\documentclass[11pt]{amsart}

\usepackage{amsfonts}

\def\ee{{\mathbf e}}
\def\hh{{\mathbf h}}

\def\Complex{{\mathbb C}}
\def\Q{{\mathbb Q}}
\def\R{{\mathbb R}}
\def\Z{{\mathbb Z}}

\def\C{{\mathcal C}}
\def\E{{\mathcal E}}

\def\ggoth{{\mathfrak g}}

\def\Card{\mathop{\mathrm{Card}}}
\def\rank{\mathop{\mathrm{rk}}}

\def\boite#1{\vrule height 15pt depth 9pt\hbox to 24pt{\hfil #1\hfil}}
\def\sboite#1{\vrule height 10pt depth 5pt\hbox to 15pt{\hfil #1\hfil}}
\def\hboite#1{\hbox to 15pt{\hfil #1\hfil}}

\newtheorem{Theorem}{Theorem}[section]
\newtheorem{Proposition}[Theorem]{Proposition}
\newtheorem{Lemma}[Theorem]{Lemma}
\newtheorem{Remarks}[Theorem]{Remark}
\newtheorem{Corollary}[Theorem]{Corollary}
\newtheorem{Example}[Theorem]{Example}
\newtheorem{Definition}[Theorem]{Definition}

\def\Item[#1]{\par\indent\hskip\parindent\llap{#1\enspace}\ignorespaces}
\def\ItemItem[#1]{\par\indent\hskip2\parindent\llap{#1\enspace}\ignorespaces}
\newenvironment{Itemize}{\vskip0.5em\begingroup}{\endgroup\vskip0.5em}

\begin{document}

\title{Centralizers of distinguished nilpotent pairs and related problems}
\author{Rupert W.T. Yu}
\address{UMR 6086 du C.N.R.S.\\
D\'epartement de Math\'ematiques\\
Universit\'e de Poitiers\\
T\'el\'eport 2 -- BP 30179\\
Boulevard Marie et Pierre Curie\\
86962 Futuroscope Chasseneuil Cedex\\
FRANCE.}
\email{yuyu@mathlabo.univ-poitiers.fr}


\begin{abstract}
In this paper, by establishing an explicit and combinatorial
description
of the centralizer of a distinguished nilpotent pair in a classical
simple Lie algebra, we solve in the classical case Panyushev's 
Conjecture which says that distinguished nilpotent pairs are wonderful,
and the classification problem on almost principal nilpotent pairs.
More precisely, we show that disinguished nilpotent pairs are wonderful
in types A, B and C, but they are not always wonderful in type D.
Also, as the corollary of the classification of almost principal 
nilpotent pairs,
we have that almost principal nilpotent pairs do not exist 
in the simply-laced case and that
the centralizer of an almost principal nilpotent
pair in a classical simple Lie algebra is always abelian.
\end{abstract}

\maketitle

\section*{Introduction}

The study of nilpotent pairs in semisimple Lie algebras is due to V. Ginzburg
\cite{[G]}. Let $\ggoth$ be a semisimple Lie algebra over an 
algebraically closed
field of characteristic zero and let $G$ be its adjoint group. Then a pair
of commuting nilpotent elements $\ee=(e_{1},e_{2})$ in $\ggoth$
is called a nilpotent
pair if there exists a pair of commuting
semisimple elements $\hh=(h_{1},h_{2})$ having
rational eigenvalues in the adjoint action such that 
$[h_{i},e_{j}]=\delta_{i,j}e_{j}$. The pair $\hh$ is called an associated
semisimple pair (or a characteristic of $\ee$ by D. Panyushev \cite{[P1]}).
The theory of nilpotent pairs can be viewed as a double counterpart of the 
theory of nilpotent orbits and Ginzburg showed in \cite{[G]} that the theory of 
principal nilpotent pairs, {\it i.e.} the simultaneous centralizer of
$\ee$ has dimension $\rank \ggoth$, has a refinement of Kostant's 
results on regular nilpotent elements in $\ggoth$. 

Although the number of $G$-orbits of nilpotent pairs is infinite, it was
shown in \cite{[G]} that, as in the classical case, the number of $G$-orbits of
principal nilpotent pairs is finite. Generalizing the notion of distinguished 
nilpotent element, Ginzburg also introduced the notion of distinguished
nilpotent pairs which are defined to be those such that the simultaneous
centralizer does not contain any semisimple elements and the
simultaneous centralizer of an associated semisimple pair is a Cartan
subalgebra, and he conjectured
that the number of $G$-orbits of distinguished nilpotent pairs is finite.
In \cite{[EP]}, A. Elashvili and D. Panyushev gave a positive answer to this
conjecture in the case of classical simple Lie algebras by giving
an explicit classification of distinguished nilpotent pairs in these
Lie algebras. 

On the other hand, Panyushev \cite{[P1]} introduced the notion of
wonderful nilpotent pairs (see Section 5 for the definition) and
he gave a proof of the finiteness of the number of $G$-orbits
of wonderful nilpotent pairs. In the pursuit of a classification-free
approach to the finiteness of the number of $G$-orbits of distinguished
nilpotent pairs, he conjectured that distinguished nilpotent pairs are
wonderful. 

In the first part of this paper, we give a complete
description of the centralizer of a distinguished
nilpotent pair for a classical simple Lie algebra $\ggoth$. This
is done by generalizing the approach of Ginzburg \cite{[G]} in the case of
${\mathrm{sl}}(V)$, using skew diagrams. This description is very explicit
and our starting point is the classification of Elashvili and Panyushev.
 
The combinatorial tool introduced here allows us to study in a purely 
diagram-theoretic manner the dimension of the centralizer of a distinguished 
nilpotent pair. In particular, we prove that the dimension of its positive
quadrant is greater than or equal to the rank of $\ggoth$, and 
we determine exactly when equality holds. An application of our results
shows that a distinguished nilpotent pair is wonderful in type A, B or C, 
but it is only so in type D under certain (necessary and sufficient)
conditions. Thus we answer the above 
conjecture completely in the classical case. 

In \cite{[P2]}, Panyushev studied another class of nilpotent pairs called
almost principal. A nilpotent pairs is called almost principal if
its centralizer has dimension $\rank \ggoth+1$. He proved that
an almost principal nilpotent pair is distinguished and wonderful and 
he gave examples of such nilpotent pairs in the non simply-laced cases. 
He suspected
that almost nilpotent pairs do not exist in the simply-laced cases
and suggested that the centralizer of an almost principal nilpotent
pair is abelian.

Using our results on the centralizers of distinguished nilpotent
pairs, we give a classification of almost principal nilpotent pairs.
As corollaries, we confirm that there are indeed no almost principal
nilpotent pairs in the simply-laced cases and also that the centralizer
of an almost principal nilpotent is effectively abelian. 

In Section 1, we define our combinatorial tool and recall the
classification of Elashvili and Panyushev. In Section 2, we prove
that the centralizer of a distinguished nilpotent pair has a basis
indexed by a certain set of pairs of skew diagrams. We study in Section 3 and 4
the combinatorics of this set of pairs of skew diagrams. In Section 5, we
give necessary and sufficient conditions for a distinguished nilpotent
pair to be wonderful. Sections 6 and 7 are devoted to the classification
of almost principal nilpotent pairs and to their centralizer. 
\vskip1em

\section*{Acknowledgment} 
The author would like to thank
Dmitri Panyushev for many discussions, Marc van Leeuwen for discussions 
on the proof of Theorem 3.2 and Pierre Torasso for his comments on 
earlier versions of this paper.

\section{Generalities on distinguished nilpotent pairs}

Let $\ee =(e_{1},e_{2})$ be a nilpotent pair, and $\hh =(h_{1},h_{2})$
an associated semisimple pair. The nilpotent pair $\ee$ is called 
{\bf distinguished} if its centralizer 
$Z (\ee ):=Z_\ggoth(e_{1})\cap Z_\ggoth(e_{2})$ in 
$\ggoth$ contains no semisimple element and the centralizer 
$Z (\hh ):=Z_\ggoth(h_{1})\cap Z_\ggoth(h_{2})$ is a Cartan 
subalgebra of $\ggoth$.

Before recalling the classification of distinguished nilpotent pairs
in classical simple Lie algebras, we shall start by defining the
combinatorial objects involved.   
 
\begin{Definition}
By a {\bf diagram}, we mean a subset $\Gamma$ of $\R^{2}$ such that
there exist $x,y\in\R$ such that $\Gamma\subset\Z^{2}+(x,y)$.
It is called {\bf connected} if given $(i,j),(k,l)\in\Gamma$, there
exists a sequence of elements $(a_{1},b_{1}), \cdots ,(a_{n},b_{n})$ in
$\Gamma$ such that

\begin{Itemize}
\Item[i)] $(a_{1},b_{1})=(i,j)$ and $(a_{n},b_{n})=(k,l)$;
\Item[ii)] for $1\leq t\leq n-1$, $(a_{t+1}-a_{t})^{2}+(b_{t+1}-b_{t})^{2}=1$.
\end{Itemize}

A {\bf skew diagram} is a connected diagram $\Gamma$ of $\R^{2}$ such that 
\vskip0.5em

\begin{Itemize}
\Item[1.] $\Gamma$ is finite;
\Item[2.] if $(i,j)\in \Gamma$ and $(i+1,j+1)\in\Gamma$, then $(i+1,j)$ and 
$(i,j+1)\in\Gamma$.
\end{Itemize}

A subset $\Gamma'$ of a skew diagram $\Gamma$ is called a {\bf skew subdiagram}
if $\Gamma'$ is a skew diagram and if $(i,j)\in\Gamma'$, then
$(i-1,j),(i,j-1)\not\in \Gamma\setminus\Gamma'$. 
\end{Definition}

\begin{Remarks}
Our definition corresponds to the usual definition of connected
skew diagrams. 
\end{Remarks}

Let $\sigma :\R^{2}\rightarrow \R^{2}$ be the central symmetric
bijection sending $(i,j)$ to $(-i,-j)$.

\begin{Definition}
A skew diagram $\Gamma$ is {\bf centrally symmetric} if $\Gamma =\sigma 
(\Gamma )$.
A subset $\Gamma'$ of a skew diagram $\Gamma$ is called a {\bf $\sigma$-skew 
subdiagram} if $\sigma (\Gamma')$ is a skew subdiagram of $\sigma (\Gamma )$. 

Note that there are three types of centrally symmetric skew diagrams:
\vskip0.5em

{\parindent = 2cm
\begin{Itemize}
\Item[{\bf Integral:}] $\Gamma\subset \Z^{2}$;
\Item[{\bf Semi-integral:}] $\Gamma\subset \Z^{2}+(1/2,0)$ or $(0,1/2)$;
\Item[{\bf Non-Integral:}] $\Gamma\subset \Z^{2}+(1/2,1/2)$;
\end{Itemize}
}
\end{Definition}

The notion of a $\sigma$-skew subdiagram will be used to describe
the centralizer of $\ee$. It corresponds to what Ginzburg called
an {\it out-subset} in \cite{[G]}, while a skew subdiagram corresponds to 
an {\it in-subset}.

%
% Basic example 1
% 
\begin{Example}
We shall represent a skew diagram by boxes. For example, 
the following is a skew diagram
which is centrally symmetric if $(i,j)=(0,0)$.
$$
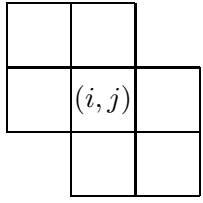
\vbox{\offinterlineskip
\halign{#&#&#&#\cr
\omit\hrulefill&\omit\hrulefill&&\cr
\boite{}&\boite{}&\boite{}&\cr
\omit\hrulefill&\omit\hrulefill&\omit\hrulefill\cr
\boite{}&\boite{$(i,j)$}&\boite{}&\boite{}\cr
\omit\hrulefill&\omit\hrulefill&\omit\hrulefill\cr
&\boite{}&\boite{}&\boite{}\cr
&\omit\hrulefill&\omit\hrulefill&\cr
}
}
$$
The subset $\{ (i-1,j),(i,j),(i,j-1)\}$ is a skew subdiagram
and $\{ (i+1,j)$, $(i,j)$, $(i,j+1)\}$ is a $\sigma$-skew subdiagram.
\end{Example}

Let us now recall the classification of distinguished nilpotent pairs in classical
simple Lie algebras of A. Elashvili and D. Panyushev.

\begin{Theorem} {\rm \cite{[EP]}}
The set of conjugacy classes of distinguished nilpotent pairs in $\ggoth$
is in bijection with:
\vskip0.5em

\noindent{\bf Type A} - 
the set of diagrams $\Gamma$ of cardinal $\rank \ggoth+1$
such that
\begin{Itemize} 
\Item[$\bullet$] $\Gamma$ is a skew diagram;
\Item[$\bullet$] the barycentre of $\Gamma$ is $(0,0)$, {\it i.e.}
$\sum_{(i,j)\in\Gamma} i=0$ and $\sum_{(i,j)\in\Gamma} j=0$
\end{Itemize}

\noindent{\bf Type B} - 
the set of pairs of centrally symmetric skew diagrams $(\Gamma^{1},\Gamma^{2})$ 
(we allow empty diagrams) such that
\begin{Itemize}
\Item[$\bullet$] $\Card \Gamma^{1}+\Card\Gamma^{2}=2\rank \ggoth +1$;  
\Item[$\bullet$] $\Gamma^{1}$ is integral and $\Gamma^{2}$ is non integral.
\end{Itemize}

\noindent{\bf Type C} - 
the set of pairs of centrally symmetric skew diagrams $(\Gamma^{1},\Gamma^{2})$ 
(we allow empty diagrams) such that
\begin{Itemize}
\Item[$\bullet$] $\Card \Gamma^{1}+\Card\Gamma^{2}=2\rank \ggoth$;  
\Item[$\bullet$] $\Gamma^{1}$ and $\Gamma^{2}$ are semi-integral
and $\Gamma^{1}\subset \Z^{2}+(0,1/2)$, $\Gamma^{2}\subset \Z^{2}+(1/2,0)$.
\end{Itemize}

\noindent{\bf Type D} - 
the set of triples $(\Gamma^{1},\Gamma^{2},\Gamma^{3})$
of centrally symmetric skew diagrams and $\epsilon\in \{ 1,2\}$ 
(we allow empty diagrams) such that
\begin{Itemize}
\Item[$\bullet$] $\Card \Gamma^{1}+\Card\Gamma^{2}+\Card \Gamma^{3}=
2\rank \ggoth$;  
\Item[$\bullet$] $\Gamma^{1}$ is non integral; 
\Item[$\bullet$] either $\Gamma^{2}$ and $\Gamma^{3}$ are empty or
they are both non empty, integral with
$\Gamma^{2}\cap\Gamma^{3}=\{ (0,0)\}\neq \Gamma^{3}$ and $\epsilon =1$.
\end{Itemize}
\end{Theorem}

We shall finish this section by expliciting the corresponding
distinguished nilpotent pairs listed in Theorem 1.5. 
Let $(\Gamma^{1},\Gamma^{2},\Gamma^{3})$ be as above, and 
$$
V=\bigoplus_{(i,j)\in\Gamma^{1}} \Complex v_{i,j}^{1}\oplus
\bigoplus_{(i,j)\in\Gamma^{2}} \Complex v_{i,j}^{2}\oplus
\bigoplus_{(i,j)\in\Gamma^{3}} \Complex v_{i,j}^{3}
$$
Define $h_{1},h_{2}\in {\mathrm{gl}}(V)$ by $h_{1}v_{i,j}^{k}=iv_{i,j}^{k}$
and $h_{2}v_{i,j}^{k}=jv_{i,j}^{k}$. 

\noindent{\bf Type A}, $\Gamma^{2}=\Gamma^{3}=\emptyset$, and
define $e_{1}v_{i,j}^{1}=v_{i+1,j}^{1}$, $e_{2}v_{i,j}^{1}=v_{i,j+1}^{1}$.
Then $\ee =(e_{1},e_{2})$ is the corresponding distinguished nilpotent
pair with $\hh =(h_{1},h_{2})$ as the associated semi-simple pair. 

\noindent{\bf Type B}, $\Gamma^{3}=\emptyset$, and
define the non-degenerate symmetric bilinear form on $V$ by
$$
(v_{i,j}^{n},v_{k,l}^{m})=\delta_{i,-k}\delta_{j,-l}\delta_{n,m}.
$$
Set 
$$
e_{1}v_{i,j}^{k}=(-1)^{i+j}v_{i+1,j}^{k}, 
e_{2}v_{i,j}^{k}=(-1)^{i+j+1}v_{i,j+1}^{k}.
$$
Then $\ee =(e_{1},e_{2})$ is the corresponding distinguished nilpotent
pair with $\hh =(h_{1},h_{2})$ as the associated semi-simple pair. 

\noindent{\bf Type C}, $\Gamma^{3}=\emptyset$, and
define the non-degenerate alternating bilinear form on $V$ by
$$
(v_{i,j}^{n},v_{k,l}^{m})=(-1)^{i+j+1/2}\delta_{i,-k}\delta_{j,-l}
\delta_{n,m}.
$$
Set 
$$
e_{1}v_{i,j}^{k}=v_{i+1,j}^{k}, 
e_{2}v_{i,j}^{k}=v_{i,j+1}^{k}.
$$
Then $\ee =(e_{1},e_{2})$ is the corresponding distinguished nilpotent
pair with $\hh =(h_{1},h_{2})$ as the associated semi-simple pair. 

\noindent{\bf Type D}, 
define the non-degenerate symmetric bilinear form on $V$ by
$$
(v_{i,j}^{n},v_{k,l}^{m})=\delta_{i,-k}\delta_{j,-l}
\delta_{n,m}.
$$
Set 
$$
e_{1}v_{i,j}^{k}=(-1)^{i+j}v_{i+1,j}^{k}, 
e_{2}v_{i,j}^{k}=(-1)^{i+j+1}v_{i,j+1}^{k}.
$$
Then $\ee =(e_{1},e_{2})$ is the corresponding distinguished nilpotent
pair with $\hh =(h_{1},h_{2})$ as the associated semi-simple pair. 

\begin{Remarks}
1. The construction for type C differs from the one in \cite{[EP]}. It seems that
our construction is more adapted for our purpose since the definitions
of the $e_{i}$'s are much simpler.

2. For type D, the parameter $\epsilon$ is not present here in the description.
This is however not important concerning the description of the centralizer.
The fact is that in this special case, there are two classes of
distinguished nilpotent pairs admitting such a description (the
corresponding bases for $V$ are not conjugate).
\end{Remarks}

\section{Description of the centralizer of a distinguished nilpotent pair}

Let $Z (\ee )$ be the centralizer of $\ee$ in $\ggoth$. 
Note that it is $h_{i}$-stable for $i=1,2$. Thus it is $\Q^{2}$-graded
and we shall denote by $Z_{p,q}(\ee )$ the homogeneous component 
of degree $(p,q)$ which is nothing but the $(h_{1},h_{2})$-eigenspace 
of eigenvalue $(p,q)$. Since $\ee$ is distinguished, the centralizer $Z(\hh )$
of $\hh$ in $\ggoth$ is a Cartan subalgebra. But $\ggoth_{0,0}=Z(\hh )$,
so $Z_{0,0}(\ee )=0$ since $Z(\ee )$ contains no semi-simple elements.

Note that by Theorem 1.5, in fact $Z_{p,q}(\ee )$ is non zero only if
$p,q$ are half integers, {\it i.e.} $2p,2q\in \Z$ and $p+q\in\Z$. 

Let us conserve the notations in Section 1 and for $1\leq k\leq 3$, set
$V_{k}=\bigoplus_{(i,j)\in \Gamma^{k}} \Complex v_{i,j}^{k}$. 

Let $x\in Z_{p,q}(\ee )$. When $k\neq l$, $V_{k}$ and $V_{l}$ 
are orthogonal with respect to the corresponding 
non-degenerate $\ggoth$-invariant bilinear form, it follows that if 
$\pi_{k}$ denotes the canonical orthogonal projection, then $x$ is the sum
$\sum_{k,l} \pi_{l}\circ x\circ\pi_{k}$. For simplicity, let us fix $k$ and $l$
and assume that $x=\pi_{l}\circ x\circ\pi_{k}+\pi_{k}\circ x\circ\pi_{l}$ if
$k\neq l$ and $x=\pi_{k}\circ x\circ\pi_{k}$ if $k=l$;
{\it i.e.} $xV_{k}\subset V_{l}$, $xV_{l}\subset V_{k}$ and 
$xV_{j}=0$ if $j\neq k,l$.

\begin{Definition}
We define $\Gamma_{x}^{k}=\{ (i,j)\in \Gamma^{k}$ such that 
$x.v_{i,j}^{k}\neq 0 \}$.
\end{Definition}

We are interested in the connected components of $\Gamma_{x}^{k}$.

\begin{Lemma}
Let $C$ be a connected component of $\Gamma_{x}^{k}$.
Then,
\begin{Itemize}
\Item[a)] $C$ is a skew subdiagram of $\Gamma^{k}$.
\Item[b)] $C+(p,q)$ is a $\sigma$-skew subdiagram of 
$\Gamma^{l}$.
\end{Itemize}
\end{Lemma}
\begin{proof}
Let $(i,j)\in C$. If $(i-1,j)$ belongs to
$\Gamma^{k}$, then $0\neq e_{1}v_{i-1,j}^{k}\in \Complex v_{i,j}^{k}$.
Hence 
$$
e_{1}xv_{i-1,j}^{k}=xe_{1}v_{i-1,j}^{k}\neq 0
$$
which implies that $(i-1,j)\in C$. The same
argument applies to $(i,j-1)$ and part a) follows.

For b), again let $(i,j)\in C$. Thus 
$0\neq xv_{i,j}^{k}\in \Complex v_{i+p,j+q}^{l}$.
If $(i+p+1,j+q)\in\Gamma^{l}$, then $e_{1}v_{i+p,j+q}^{l}\neq 0$. 
Hence 
$$
xe_{1}v_{i,j}=e_{1}xv_{i,j}\neq 0.
$$
It follows that $xv_{i+1,j}\neq 0$, $(i+1,j)\in C$
and $(i+1+p,j+q)\in C+(p,q)\subset \Gamma^{l}$.
The same argument applies to $(i+p,j+q+1)$ and we have proved
that $C+(p,q)$ is a $\sigma$-skew subdiagram of $\Gamma^{l}$.
\end{proof}

Note that $C$ and $C+(p,q)$ are of the same
shape.

Denote by $\C_{x}^{k}$ the set of connected components 
of $\Gamma_{x}^{k}$.

\begin{Proposition}
Let $\ggoth$ be of type $A$, then $\Gamma=\Gamma^{1}$. 
For $C\in \C_{x}^{1}$, define $x_{C}$ to be the obvious restriction of $x$
to $V_{C}=\bigoplus_{(i,j)\in C}\Complex v_{i,j}^{1}$. 
Then $x_{C}\in Z_{p,q}(\ee )$ and $x=\sum_{C\in\C_{x}^{1}} x_{C}$.
\end{Proposition}
\begin{proof}
Since $x$ commutes with the $e_{i}$'s, it is clear that $x_{C}$ is nilpotent
and commutes also 
with the $e_{i}$'s and so the first statement is clear. The last statement is 
also clear since we are taking connected components.
\end{proof}

For types $B$, $C$ and $D$, the $\Gamma^{j}$'s are centrally symmetric
and we have the following lemma.

\begin{Lemma}
Let $\ggoth$ be of type other than $A$, then for any $C\in\C_{x}^{k}$, 
$\sigma (C)-(p,q)\in \C_{x}^{l}$.
\end{Lemma}
\begin{proof}
By Lemma 2.2, $C+(p,q)$ is a $\sigma$-skew subdiagram of $\Gamma^{l}$.
Since $\Gamma^{l}$ is centrally symmetric, $\sigma (C)-(p,q)=\sigma (C+(p,q))$
is a skew subdiagram of $\Gamma^{l}$.

Now let $(i,j)\in C$ and denote by $(\, ,\, )$, the $\ggoth$-invariant
non-degenerate bilinear form corresponding to the type of $\ggoth$ as
explained in Section 1. We have
$$
(xv_{i,j}^{k},v_{-i-p,-j-q}^{l})+(v_{i,j}^{k},xv_{-i-p,-j-q}^{l})=0.
$$
Thus $xv_{-i-p,-j-q}^{l}\neq 0$ and $(-i-p,-j-q)\in \Gamma_{x}^{l}$.

To finish the proof, it suffices to show that $\sigma (C)-(p,q)$ is
a connected component of $\Gamma_{x}^{l}$.
Let $(i,j)\in \sigma (C)-(p,q)$ be such that $(i+1,j)\in\Gamma_{x}^{l}$.
We have $0\neq xv_{i+1,j}^{l}\in \Complex v_{i+1+p,j+q}^{k}$. By central
symmetry, this means that $(-i-1-p,-j-q)\in\Gamma^{k}$. This in turn
implies that $(-i-1-p,-j-q)\in C$ since $(-i-p,-j-q)\in C$ and $C$
is a skew subdiagram. Thus $(i+1,j)\in \sigma (C)-(p,q)$.

The same argument applies for $(i,j+1)$ and this proves that
$\sigma (C)-(p,q)$ is a connected component of $\Gamma_{x}^{l}$.
\end{proof}
 
\begin{Definition}
We define the set 
$$
\C_{x}^{k,l} :=\{ \{ C, C'\} \,\mid\, C\in \C_{x}^{k}, C'\in \C_{x}^{l} \,
{\mathrm{and}}\, \sigma (C)=C'+(p,q)\}.
$$
\end{Definition}

\begin{Remarks}
First of all, note that we take subsets instead of pairs to avoid repetitions 
when $k=l$.
Also it is clear that $\C_{x}^{k,l}$ and $\C_{x}^{l,k}$
are identical. 
\end{Remarks}

\begin{Proposition}
Let $\ggoth$ be of type other than $A$, $\{ C,C'\}\in \C_{x}^{k,l}$
and $x_{C,C'}$ the restriction of $x$ to 
$\bigoplus_{(i,j)\in C}\Complex v_{i,j}^{k} \oplus
\bigoplus_{(i,j)\in C'}\Complex v_{i,j}^{l}$. 
Then $x_{C,C'}\in Z_{p,q}(\ee )$ and $x=\sum_{\{ C,C'\}\in\C_{x}^{k,l}} 
x_{C,C'}$.
\end{Proposition}
\begin{proof}
It is clear that $x_{C,C'}$ commutes with the $e_{i}$'s and a direct
computation shows that $x_{C,C'}\in \ggoth$. The second statement
follows from Lemma 2.4 (which is also true if we exchange $k$ et $l$)
since we are considering connected components.
\end{proof}

Note that $x_{C,C'}=x_{C',C}$.

It is now clear how to define an indexing set for a set of
generators of $Z_{p,q}(\ee )$.

\begin{Definition} 
Let $\E_{p,q}(\ee)^{k}$ be the set of skew subdiagrams $C$ in 
$\Gamma^{k}$ such that $C'=C+(p,q)$ is a $\sigma$-skew subdiagram of 
$\Gamma^{k}$.

For $k\neq l$, 
let $\E_{p,q}(\ee)^{k,l}$ be the set of pairs of skew subdiagrams 
$(C_{1},C_{2})$ such that $C_{1}$ is a skew subdiagram of $\Gamma^{k}$,
$C_{2}$ is a skew subdiagram of $\Gamma^{l}$ and  
$C_{1}+(p,q)=\sigma (C_{2})$.

For $k=l$, $p,q$ are integers. 
If $p+q$ is odd, then let $\E_{p,q}(\ee)^{k,k}$ be the set of subsets of skew 
subdiagrams  $\{ C_{1},C_{2}\}$ in $\Gamma^{k}$ 
such that $C_{1}+(p,q)=\sigma (C_{2})$.
If $p+q$ even, then let $\E_{p,q}(\ee)^{k,k}$ be the set of subsets  
of skew subdiagrams $\{ C_{1},C_{2}\}$ in $\Gamma^{k}$
such that $C_{1}\neq C_{2}$ and $C_{1}+(p,q)=\sigma (C_{2})$.

Set $\E_{p,q}(\ee)$ to be $\E_{p,q}(\ee )^{1}$ if $\ggoth$ is of type
$A$ with $\Gamma=\Gamma^{1}$
and to be the disjoint union of the $\E_{p,q}(\ee )^{k,l}$'s with $k\leq l$
if $\ggoth$ is of type $B$, $C$ or $D$. 
\end{Definition}

\begin{Theorem}
For $(p,q)\neq (0,0)$, the subspace $Z_{p,q}(\ee)$ has a basis 
indexed by the set $\E_{p,q}(\ee )$.
\end{Theorem}
\begin{proof}
We shall first show that if $x\in \E_{p,q}(\ee )$ and
$C_{1}\in \C_{x}^{1}$ (resp. $\{ C_{1},C_{2}\}\in \C_{x}^{k,l}$), 
then the element $x_{C_{1}}$ (resp. $x_{C_{1},C_{2}}$) 
as defined in Propositions 2.3 (resp. Proposition 2.7)
is uniquely determined by $C_{1}$ (resp. $C_{1},C_{2}$) 
up to a scalar.

Let $e_{1}v_{i,j}^{k}=a_{i,j}^{k}v_{i+1,j}^{k}$,
$e_{2}v_{i,j}^{k}=b_{i,j}^{k}v_{i+1,j}^{k}$ and if appropriate, we let
$(v_{i,j}^{k}$, $v_{-i,-j}^{k})$ $=d_{i,j}^{k}$.

Let $x$ be in $Z_{p,q}(\ee )$ with the appropriate property.
For $(i,j)\in \Gamma_{x}^{k}$, $xv_{i,j}^{k}=\lambda_{i,j}^{k}v_{i+p,j+q}^{l}$. 

Now let us fix $k$ and $l$ and let $(i,j)\in C_{1}\in\C_{x}^{k,l}$.
If $(i+1,j)\in C_{1}$, then 
$e_{1}x_{C_{1}}v_{i,j}^{k}=x_{C_{1}}e_{1}v_{i,j}^{k}$
implies that 
$$
a_{i+p,j+q}^{l}\lambda_{i,j}^{k}=\lambda_{i+1,j}^{k}a_{i,j}^{k}.
$$
Similarly, if $(i,j+1)\in C_{1}$, then
$$
b_{i+p,j+q}^{l}\lambda_{i,j}^{k}=\lambda_{i,j+1}^{k}b_{i,j}^{k}.
$$
In the case of type A and C, $a_{i,j}^{k}=b_{i,j}^{k}=1$ and we have
$$
\lambda_{i+1,j}^{k}=\lambda_{i,j}^{k}=\lambda_{i,j+1}^{k}.
$$
In the case of type B and D, $a_{i,j}^{k}=(-1)^{i+j}=-b_{i,j}^{k}$
and we have
$$
\lambda_{i+1,j}^{k}=(-1)^{p+q}\lambda_{i,j}^{k}=\lambda_{i,j+1}^{k}.
$$
Also, we have for types B, C and D,
$$
(xv_{i,j}^{k},v_{-i-p,-j-q}^{l})+(v_{i,j}^{k},xv_{-i-p,-j-q}^{l})=0
$$
which implies that
$$
\lambda_{i,j}^{k}d_{i+p,j+q}^{l}+\lambda_{-i-p,-j-q}^{l}d_{i,j}^{k}=0.
$$
When $k=l$, we have $p,q\in\Z$.
If $p+q$ is even, then we must have 
$C_{1}\neq C_{2}$ because from the classification of Section 1,
$d_{i,j}^{k}=1$ for types B, D and $d_{i,j}^{k}=(-1)^{i+j+1/2}$ for
type C (one should treat the case where $C_{1}=\{(i,j)\}$ separately). 

Since $p,q$ are fixed and $C_{1}$ is connected, we observe that the restriction
of $x$ on $V_{C_{1}}$ and $V_{C_{2}}$ are uniquely determined by the value of 
$\lambda_{i,j}^{k}$ ($i,j,k$ are fixed here). It follows that
$x_{C_{1}}$ (resp. $x_{C_{1},C_{2}}$) is uniquely determined by $C_{1}$
(resp. $(C_{1},C_{2})$) up to a scalar. 

By Propositions 2.3 and 2.7, we have that $Z_{p,q}(\ee )$ is spanned
by a set of elements indexed by a subset of $\E_{p,q}(\ee )$.
To finish the proof of this theorem, il suffices to construct an element of
$Z_{p,q}(\ee )$ for each element of $\E_{p,q}(\ee )$ and show that
they form a basis of $Z_{p,q}(\ee )$. Recall that
$Z_{p,q}(\ee )$ is non zero only if $2p,2q\in\Z$ and $p+q\in\Z$.

\noindent{\bf Type A.} Let $\Gamma$ be a skew diagram with barycentre $(0,0)$
and $C\in\E_{p,q}(\Gamma)$.
Define $x_{C}\in {\mathrm{gl}}(V)$ by
$$
x_{C}v_{i,j}=\left\{
\begin{array}{cc}
v_{i+p,j+q}&{\mathrm{if}}\, (i,j)\in C\\
0&{\mathrm{otherwise.}}
\end{array}\right.
$$
Then a simple verification shows that $x_{C}$ commutes with
$e_{1}$ and $e_{2}$. Further $x_{C}$ is nilpotent since $(p,q)\neq (0,0)$
and therefore $x_{C}\in {\mathrm{sl}}(V)$.

\noindent{\bf Types B and D.} Let $\Gamma^{k}$, $\Gamma^{l}$ be
centrally symmetric skew diagrams. Let 
$\{ C,C'\}\in\E_{p,q}(\Gamma^{k},\Gamma^{l})$ be
such that $k\neq l$ or $C \neq C'$. 
Define $x_{C,C'}\in {\mathrm{gl}}(V)$ by
$$
x_{C,C'}v_{i,j}^{m}=\left\{
\begin{array}{cc}
v_{i+p,j+q}^{l}&{\mathrm{if}}\, (i,j)\in C, m=k\\
- v_{i+p,j+q}^{k}&{\mathrm{if}}\, (i,j)\in C', m=l\\
0&{\mathrm{otherwise}}
\end{array}\right.
$$
if $p+q$ is even; and if $p+q$ is odd, we set
$$
x_{C,C'}v_{i,j}^{m}=\left\{
\begin{array}{cc}
(-1)^{i+j}v_{i+p,j+q}^{l}&{\mathrm{if}}\, (i,j)\in C, m=k\\
(-1)^{i+j}v_{i+p,j+q}^{k}&{\mathrm{if}}\, (i,j)\in C', m=l\\
0&{\mathrm{otherwise.}}
\end{array}\right.
$$
Again a simple verification shows that $x_{C,C'}$ commutes with
$e_{1}$ and $e_{2}$, and that $x_{C,C'}$ is nilpotent since $(p,q)\neq (0,0)$.
Finally, we verify easily that for $(i,j)\in C$,
$$
(x_{C,C'}v_{i,j}^{k},v_{-i-p,-j-q}^{l})+
(v_{i,j}^{k},x_{C,C'}v_{-i-p,-j-q}^{l})=0.
$$
Therefore $x_{C,C'}\in Z_{p,q}(\ee )$.

Now suppose that $l=k$ and $p+q$ is odd. Let $C\in\E_{p,q}(\Gamma^{k})$ be
such that $C=\sigma(C)-(p,q)$.
Define $x_{C,C}\in {\mathrm{gl}}(V)$ by
$$
x_{C,C}v_{i,j}^{m}=\left\{
\begin{array}{cc}
(-1)^{i+j}v_{i+p,j+q}^{k}&{\mathrm{if}}\, (i,j)\in C, m=k\\
0&{\mathrm{otherwise.}} 
\end{array}\right.
$$
The same verification works and therefore $x_{C,C}\in Z_{p,q}(\ee )$.

\noindent{\bf Type C.} Let $\Gamma^{k}$, $\Gamma^{l}$ be
centrally symmetric skew diagrams. Let 
$\{ C,C'\}\in\E_{p,q}(\Gamma^{k},\Gamma^{l})$ be
such that $k\neq l$ or $C\neq C'$. 
Define $x_{C,C'}\in {\mathrm{gl}}(V)$ by
$$
x_{C,C'}v_{i,j}^{m}=\left\{
\begin{array}{cc}
v_{i+p,j+q}^{l}&{\mathrm{if}}\, (i,j)\in C, m=k\\
(-1)^{p+q+1}v_{i+p,j+q}^{k}&{\mathrm{if}}\, (i,j)\in C', m=l\\
0&{\mathrm{otherwise.}} 
\end{array}\right.
$$
Again a simple verification shows that $x_{C,C'}$ commutes with
$e_{1}$ and $e_{2}$, and that $x_{C,C'}$ is nilpotent since $(p,q)\neq (0,0)$.
Finally, we verify easily that for $(i,j)\in C$,
$$
(x_{C,C'}v_{i,j}^{k},v_{-i-p,-j-q}^{l})+
(v_{i,j}^{k},x_{C,C'}v_{-i-p,-j-q}^{l})=0.
$$
Therefore $x_{C,C'}\in Z_{p,q}(\ee )$.

Now suppose that $l=k$ and $p+q$ is odd. Let $C\in\E_{p,q}(\Gamma^{k})$ be
such that $C=\sigma(C)-(p,q)$.
Define $x_{C,C}\in {\mathrm{gl}}(V)$ by
$$
x_{C,C}v_{i,j}^{m}=\left\{
\begin{array}{cc}
v_{i+p,j+q}^{l}&{\mathrm{if}}\, (i,j)\in C, m=k\\
0&{\mathrm{otherwise.}}  
\end{array}\right.
$$
The same verification works and therefore $x_{C,C}\in Z_{p,q}(\ee )$.
\vskip1em

Thus we have constructed a set of elements of $Z_{p,q}(\ee )$ indexed
by the set $\E_{p,q}(\ee )$. Finally,
since they are defined on connected components or pairs of 
connected components, they are clearly linearly independant. 
The proof is now complete.
\end{proof}

%
% Example 
%
\begin{Example}
Let us end this section by expliciting the centralizer of the
distinguished nilpotent pair in $B_{3}$ of the centrally symmetric skew diagram
$\Gamma$ in Example 1.4. It has 
in fact a basis given by:
$$
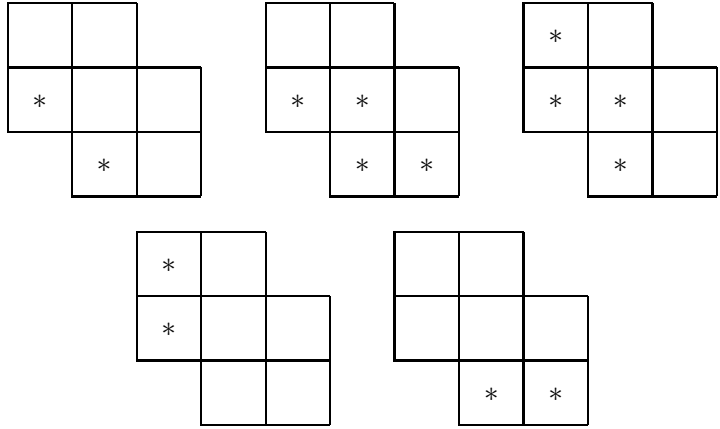
\vbox{\offinterlineskip
\halign{#&#&#&#&#&#&#&#&#&#&#&#&#\cr
\omit\hrulefill&\omit\hrulefill&&&
\omit\hrulefill&\omit\hrulefill&&&
\omit\hrulefill&\omit\hrulefill&&&\cr
\boite{}&\boite{}&\boite{}&&
\boite{}&\boite{}&\boite{}&&
\boite{$*$}&\boite{}&\boite{}&&\cr
\omit\hrulefill&\omit\hrulefill&\omit\hrulefill&&
\omit\hrulefill&\omit\hrulefill&\omit\hrulefill&&
\omit\hrulefill&\omit\hrulefill&\omit\hrulefill&&\cr
\boite{$*$}&\boite{}&\boite{}&\boite{}&
\boite{$*$}&\boite{$*$}&\boite{}&\boite{}&
\boite{$*$}&\boite{$*$}&\boite{}&\boite{}&\cr
\omit\hrulefill&\omit\hrulefill&\omit\hrulefill&&
\omit\hrulefill&\omit\hrulefill&\omit\hrulefill&&
\omit\hrulefill&\omit\hrulefill&\omit\hrulefill&&\cr
&\boite{$*$}&\boite{}&\boite{}&
&\boite{$*$}&\boite{$*$}&\boite{}&
&\boite{$*$}&\boite{}&\boite{}&\cr
&\omit\hrulefill&\omit\hrulefill&&
&\omit\hrulefill&\omit\hrulefill&&
&\omit\hrulefill&\omit\hrulefill&&\cr
}
}
$$
$$
\vbox{\offinterlineskip
\halign{#&#&#&#&#&#&#&#&#&#&#&#&#\cr
\omit\hrulefill&\omit\hrulefill&&&
\omit\hrulefill&\omit\hrulefill&&&\cr
\boite{$*$}&\boite{}&\boite{}&&
\boite{}&\boite{}&\boite{}&&\cr
\omit\hrulefill&\omit\hrulefill&\omit\hrulefill&&
\omit\hrulefill&\omit\hrulefill&\omit\hrulefill&&\cr
\boite{$*$}&\boite{}&\boite{}&\boite{}&
\boite{}&\boite{}&\boite{}&\boite{}&\cr
\omit\hrulefill&\omit\hrulefill&\omit\hrulefill&&
\omit\hrulefill&\omit\hrulefill&\omit\hrulefill&&\cr
&\boite{}&\boite{}&\boite{}&
&\boite{$*$}&\boite{$*$}&\boite{}&\cr
&\omit\hrulefill&\omit\hrulefill&&
&\omit\hrulefill&\omit\hrulefill&&\cr
}
}
$$
where the ``star'' boxes denote the corresponding pairs of
skew subdiagrams. Let $x_{1}$, $x_{2}$, $x_{3}$, $x_{4}$, $x_{5}$ be 
the corresponding
elements in ${\mathrm{so}}_{7}$ (from left to right and top to bottom) defined by
$x_{1}v_{0,-1}=v_{1,0}$, $x_{1}v_{-1,0}=-v_{0,1}$ and zero otherwise;
$x_{2}v_{i,j}=(-1)^{i+j+1}v_{i,j+1}$, $x_{3}v_{i,j}=(-1)^{i+j+1}v_{i+1,j}$,
$x_{4}v_{i,j}=(-1)^{i+j+1}v_{i+2,j-1}$, $x_{5}v_{i,j}=(-1)^{i+j+1}v_{i-1,j+2}$ 
for all $(i,j)$ with $v_{i,j}=0$ if $(i,j)\not\in \Gamma$. A simple 
verification shows that this is a basis of $Z (\ee )$.
\end{Example}

\section{Description of the positive quadrant via skew diagrams}

Since the sets $\E_{p,q}(\ee )^{k}$ and $\E_{p,q}(\ee )^{k,l}$ depend only on 
the connected skew diagrams $\Gamma^{k}$ and $\Gamma^{l}$, we can study it in 
a combinatorial way. For a skew diagram $\Gamma$, we define 
%\vskip0.5em

\begin{Itemize}
\Item[$\bullet$] $\E_{p,q}(\Gamma )$ to be the set of skew subdiagrams $C$
in $\Gamma$ such that $C+(p,q)$ is a $\sigma$-skew subdiagram in $\Gamma$.
\end{Itemize}
%\vskip0.5em

\noindent and for $\Gamma$, $\Gamma'$ two distinct
centrally symmetric skew diagrams, we define

\begin{Itemize} 
\Item[$\bullet$] $\E_{p,q}(\Gamma,\Gamma')$ to be the 
set of pairs of skew subdiagrams $(C,C')$ such that $C$ is a skew subdiagram
of $\Gamma$, $C'$ is a skew subdiagram of $\Gamma'$ and $C+(p,q)=\sigma (C')$.

\Item[$\bullet$] $\E_{p,q}(\Gamma,\Gamma )$ to be the 
set of subsets of skew subdiagrams $\{ C,C'\}$ of $\Gamma$ such that 
$C+(p,q)=\sigma (C')$ and if $p+q$ is even, then $C\neq C'$.
\end{Itemize}

In this section, we are interested in the description of the positive 
quadrant, that is when $p,q\geq 0$. Let us fix $p,q\geq 0$. Note that
if $\Gamma$, $\Gamma'$ are both integral or non integral, then
$p,q\in\Z$.

\begin{Lemma}
Let $\Gamma$ be a skew diagram, then the cardinality of $\E_{p,q}(\Gamma )$ 
is the number of connected components of $(\Gamma -(p,q)) \cap \Gamma$.
In particular $\E_{0,0}(\Gamma )$ is of cardinal $1$.
\end{Lemma}
\begin{proof}
Since $p,q\geq 0$, it is clear that a connected component $C$ of 
$(\Gamma -(p,q)) \cap \Gamma$ is a skew subdiagram of $\Gamma$ and 
that $C+(p,q)$ is a $\sigma$-skew subdiagram of $\Gamma$.  
\end{proof}

Let $\Gamma$ be a skew diagram. An element
$(i,j)\in\Gamma$ is called a northeast (resp. southwest) corner if  
$(i+1,j)$ and $(i,j+1)$ (resp. $(i-1,j)$ and $(i,j-1)$)
are not in $\Gamma$. It is called a northeast (resp. southwest)
angle if $(i+1,j+1)$ (resp. $(i-1,j-1)$) is not in $\Gamma$ but $(i+1,j)$
and $(i,j+1)$ (resp. $(i-1,j)$ and $(i,j-1)$) are in $\Gamma$. 

%
% Basic example 2
% 
For example, let $\Gamma$ be the following skew diagram,
$$
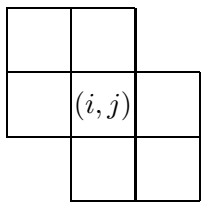
\vbox{\offinterlineskip
\halign{#&#&#&#\cr
\omit\hrulefill&\omit\hrulefill&&\cr
\boite{}&\boite{}&\boite{}&\cr
\omit\hrulefill&\omit\hrulefill&\omit\hrulefill\cr
\boite{}&\boite{$(i,j)$}&\boite{}&\boite{}\cr
\omit\hrulefill&\omit\hrulefill&\omit\hrulefill\cr
&\boite{}&\boite{}&\boite{}\cr
&\omit\hrulefill&\omit\hrulefill&\cr
}
}
$$
then $(i,j)$ is a southwest angle which is also a northeast angle,
while $(i-1,j)$ and $(i,j-1)$ are southwest corners, and 
$(i+1,j)$, $(i,j+1)$ are northeast corners.

\begin{Theorem}
Let $\Gamma$ be a skew diagram of cardinal $n$. Then
the cardinal of the set $\bigcup_{p,q\geq 0, (p,q)\neq (0,0)} 
\E_{p,q}(\Gamma )$ is equal to $n-1$.
\end{Theorem}
\begin{proof}
First observe that $\E_{0,0}(\Gamma )$ has cardinal $1$.

By the previous lemma, we are reduced to counting the number of
connected components $e_{p,q}(\Gamma )$ of $(\Gamma -(p,q))\cap \Gamma$. 
We shall proceed by induction 
on the cardinality of $\Gamma$. The result is trivial
if the cardinality of $\Gamma$ is $1$. 

Now let $n >1$ be the cardinality of $\Gamma$. Let us first
suppose that there exist a northeast corner 
$(i,j)\in\Gamma$ such that $\Gamma'=\Gamma\setminus \{ (i,j) \}$ 
is again a skew diagram. 

We shall compare the connected components of $(\Gamma'-(p,q))\cap \Gamma'$
and $(\Gamma -(p,q))\cap \Gamma$. Note that the former is obtained
from the latter by removing $(i-p,j-q)$ set-theoretically. It follows by
an easy observation that,

\begin{Itemize}
\Item[1.] $e_{p,q}(\Gamma')=e_{p,q}(\Gamma )+1$ if  $(i,j)-(p,q)$ is a
southwest angle of $\Gamma'$, which is also a southwest angle of
$\Gamma$.
\Item[2.] $e_{p,q}(\Gamma')=e_{p,q}(\Gamma )-1$ if  $(i,j)-(p,q)$ is a 
southwest corner of $\Gamma'$, which is also a southwest corner of
$\Gamma$.
\Item[3.] $e_{p,q}(\Gamma')=e_{p,q}(\Gamma )$ otherwise.
\end{Itemize}

%
% Basic example 3
% 
For example to illustrate 1, let $\Gamma$ be the following skew diagram,
$$
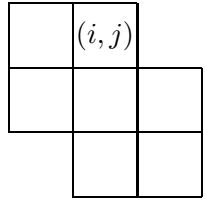
\vbox{\offinterlineskip
\halign{#&#&#&#\cr
\omit\hrulefill&\omit\hrulefill&&\cr
\boite{}&\boite{$(i,j)$}&\boite{}&\cr
\omit\hrulefill&\omit\hrulefill&\omit\hrulefill\cr
\boite{}&\boite{}&\boite{}&\boite{}\cr
\omit\hrulefill&\omit\hrulefill&\omit\hrulefill\cr
&\boite{}&\boite{}&\boite{}\cr
&\omit\hrulefill&\omit\hrulefill&\cr
}
}
$$
If we remove $(i,j)$ and take $(p,q)=(0,1)$, then $(\Gamma-(p,q))\cap
\Gamma$ looks like
$$
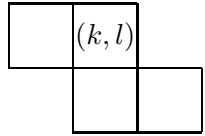
\vbox{\offinterlineskip
\halign{#&#&#&#\cr
\omit\hrulefill&\omit\hrulefill\cr
\boite{}&\boite{$(k,l)$}&\boite{}\cr
\omit\hrulefill&\omit\hrulefill&\omit\hrulefill\cr
&\boite{}&\boite{}&\boite{}\cr
&\omit\hrulefill&\omit\hrulefill&\cr
}
}
$$
where $(k,l)=(i,j-1)$ is a southwest angle of $\Gamma$. 
So when we remove this box, we recover $(\Gamma'-(p,q))\cap \Gamma'$ and
the number of connected components is increased by $1$. 

\vskip0.5em

We are therefore reduced to counting southwest corners and southwest
angles of $\Gamma$ southwest of $(i,j)$.
Let $N_{c}$ be the number of southwest corners $(k,l)$ in $\Gamma$ 
such that $k\leq i$ and $l\leq j$,
and $N_{a}$ the number of southwest angles $(k,l)$ in $\Gamma$ 
such that $k\leq i$ and $l\leq j$.
Then it follows from the above that
$$
\sum_{p,q\geq 0, (p,q)\neq (0,0)} 
e_{p,q}(\Gamma )=\sum_{p,q\geq 0, (p,q)\neq (0,0)}e_{p,q}(\Gamma')
+N_{c}-N_{a}.
$$
Now $\Gamma$ is a skew diagram, so each southwest angle is above a 
southwest corner. Further, the 
leftmost southwest corner is not below any southwest angle. We conclude
by induction that
$$
\sum_{p,q\geq 0, (p,q)\neq (0,0)} e_{p,q}(\Gamma )=
\sum_{p,q\geq 0, (p,q)\neq (0,0)}e_{p,q}(\Gamma')+1=n-1.
$$

We are left with the case where removing a northeast corner
in $\Gamma$ gives $2$ skew diagrams. This means that $\Gamma$ 
is of the form:

%
% Basic example 4
%
$$
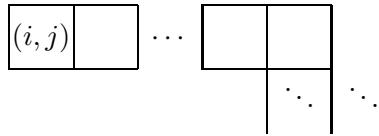
\vbox{\offinterlineskip
\halign{#&#&#&#&#&#\cr
\omit\hrulefill&\omit\hrulefill&&\omit\hrulefill&\omit\hrulefill&\cr
\boite{$(i,j)$}&\boite{}&\boite{$\cdots$}&\boite{}&\boite{}&\boite{}\cr
\omit\hrulefill&\omit\hrulefill&&\omit\hrulefill&\omit\hrulefill&\cr
&&&&\boite{$\ddots$}&\boite{$\ddots$}\cr
}
}
$$

Let $(i,j)$ be the leftmost, {\it i.e.} $i$ minimal, 
southwest corner of $\Gamma$,
and let $\Gamma'=\Gamma\setminus\{ (i,j)\}$. Then it is clear
that $e_{p,q}(\Gamma )=e_{p,q}(\Gamma' )$ if $q > 0$.
Let $(i+r,j)$ be the leftmost northeast corner of $\Gamma$.
Then 
$$
e_{r,0}(\Gamma )=\left\{
\begin{array}{cc}
e_{r,0}(\Gamma' )\hfil & {\mathrm{if} }\, p\neq r,\\
e_{r,0}(\Gamma')+1 & {\mathrm{if}}\, p=r.
\end{array}\right.
$$
It follows that $\sum_{p,q\geq 0, (p,q)\neq (0,0)}
e_{p,q}(\Gamma )=\sum_{p,q\geq 0, (p,q)\neq (0,0)}
e_{p,q}(\Gamma')+1=n-1$.
\end{proof}

Let us now turn to the study of $\E_{p,q}(\Gamma ,\Gamma )$.

\begin{Lemma}
Let $\Gamma$ be a centrally symmetric skew diagram, and let $C$ be an
element of $\E_{p,q}(\Gamma )$. Set $C'=\sigma (C)-(p,q)$. Then
$C'=C$ or $C\cap C'=\emptyset$.
\end{Lemma}
\begin{proof}
First of all, we note that $C'+(p,q)$ is a $\sigma$-skew subdiagram of 
$\Gamma$, and since $\Gamma$ is centrally symmetric, $C'$ is a skew subdiagram 
of $\Gamma$. Thus $C'\in\E_{p,q}(\Gamma )$.

Let us suppose that $C\cap C'$ is non empty and $C$ is not a subset of $C'$. 
Then there exists $(i,j)\in C\setminus C'$ such that $(i-1,j)$ or $(i,j-1)$
is in $C'$. This is not possible because $(i+p,j+q)\in \Gamma\setminus 
C'+(p,q)$, while $(i-1+p,j+q)$ or $(i+p,j-1+q)$ is in $C'+(p,q)$ which, 
as noted above, is a $\sigma$-skew subdiagram. The lemma now follows.
\end{proof}

\begin{Remarks}
It follows from the previous lemma that $\E_{p,q}(\Gamma )$ is the disjoint
union of those $C$ such that $C=C'$ and subsets $\{ C,C'\}$ with 
$C\cap C'=\emptyset$, where $C'=\sigma (C)-(p,q)$. Further, by central
symmetry, there is at most one $C$ in $\E_{p,q}(\Gamma )$ such that $C=C'$.
\end{Remarks}

\begin{Theorem}
Let $\Gamma$ be a centrally symmetric skew diagram of cardinal $n$, 
then the cardinality of $\bigcup_{p,q\geq 0, (p,q)\neq (0,0)}
\E_{p,q}(\Gamma ,\Gamma)=
\bigcup_{p,q\geq 0} \E_{p,q}(\Gamma ,\Gamma )$
is equal to $[n/2]$, where $[n/2]$ denotes the integer part of $n/2$.
\end{Theorem}
\begin{proof}
First observe that $(C,C') \in \E_{0,0}(\Gamma,\Gamma )$ implies that
$C=C'$. But $p+q$ is even, so $\E_{0,0}(\Gamma,\Gamma )=\emptyset$.

We proceed by induction on $n$ as in Theorem 3.2. For $n=1,2$, the result
is clear. Let $n > 2$ and
let $(i,j)$ be the northeast corner of $\Gamma$ such that $i$ is minimal.
Let $e_{p,q}(\Gamma,\Gamma )$ denote the cardinality of $\E_{p,q}(\Gamma,\Gamma )$.

Let us first take care of the case where $\Gamma$ is rectangular, {\it i.e.}
$\Gamma$ has exactly one southwest corner since $\Gamma$ is centrally symmetric.
Then every $C\in\E_{p,q}(\Gamma )$
is rectangular and $C=\sigma (C)-(p,q)$. So we have only to count the number
of $p,q$'s such that $p+q$ is odd. This is clearly $[n/2]$.

So let us suppose that $\Gamma$ is not rectangular.

Let $\Gamma'=\Gamma \setminus \{(i,j),(-i,-j)\}$ and let us suppose that
$\Gamma'$ is again a skew diagram. 

We shall use a similar argument as in Theorem 3.2, but since we are
dealing with centrally symmetric skew diagrams, we need to have
a detailed analysis on the induction process.

Denote by
$\Delta_{p,q}=(\Gamma - (p,q))\cap \Gamma$
and $\Delta_{p,q}'=(\Gamma' - (p,q))\cap \Gamma'$.
Note that $\Delta_{p,q}'=\Delta_{p,q}\setminus \{(-i,-j),(i-p,j-q)\}$.

Let $C$ be a connected component of $\Delta_{p,q}$, then we have the 
following 3 lemmas:

\begin{Lemma}
We have $C\cap \{(-i,-j),(i-p,j-q)\}=\emptyset$ if and only if 
$C$ is a connected component of $\Delta_{p,q}'$.
\end{Lemma}

\begin{Lemma}
Let $C\cap \{(-i,-j),(i-p,j-q)\}=\{(i-p,j-q)\}$ and $(-i,-j)\neq
(i-p,j-q)$, then we have $C\neq \sigma (C)-(p,q)$, 
$(\sigma (C)-(p,q)) \cap \{ (-i,-j), (i-p,j-q)\} = \{ (-i,-j)\}$ and
\begin{Itemize}
\Item[a)] $C\setminus \{(i-p,j-q)\}=\emptyset$ if and only if 
$(i-p,j-q)$ is a southwest corner of $\Gamma$ southwest of $(i,j)$.
\Item[b)] $C\setminus \{(i-p,j-q)\} = C'$ is a connected component of 
$\Delta_{p,q}'$ if and only if
$(i-p,j-q)$ is neither a southwest corner nor
a southwest angle of $\Gamma$ southwest of $(i,j)$. In particular, 
$\sigma (C')-(p,q)\neq C'$.
\Item[c)] $C\setminus \{(i-p,j-q)\} = C'\cup C''$ is the disjoint 
union of two distinct connected components of $\Delta_{p,q}'$ if and only if
$(i-p,j-q)$ is a southwest angle of $\Gamma$ southwest of $(i,j)$. 
In particular, $\sigma (C')-(p,q)\neq C'$,
$\sigma (C'')-(p,q)\neq C''$ and $\sigma (C')-(p,q)\neq C''$.
\end{Itemize}
\end{Lemma}

\begin{Lemma}
Let $C\cap \{(-i,-j),(i-p,j-q)\}=\{(-i,-j),(i-p,j-q)\}$, then
$C=\sigma (C)-(p,q)$ and
\begin{Itemize}
\Item[a)] $C\setminus \{(-i,-j),(i-p,j-q)\}=\emptyset$ if and only if either
$(-i,-j)\in \{ (i-p,j-q),(i-p-1,j-q)\}$ or $(-i,-j)=(i-p,j-q-1)$ and 
$(i-p,j-q)$ is not a southwest angle of $\Gamma$ southwest of $(i,j)$.
\Item[b)] $C\setminus \{(-i,-j),(i-p,j-q)\}=C'\cup C''$ is the disjoint 
union of two distinct connected components of $\Delta_{p,q}'$ such
that $C'=\sigma (C'')-(p,q)$ if and only if either $(-i,-j)=(i-p-1,j-q-1)$ or
$(-i,-j)=(i-p,j-q-1)$ and $(i-p,j-q)$ is a southwest angle of $\Gamma$ 
southwest of $(i,j)$.
\Item[c)] $C\setminus \{(-i,-j),(i-p,j-q)\}=C'\cup C''\cup C'''$ is the 
disjoint 
union of three distinct connected components of $\Delta_{p,q}'$ such that
$C'=\sigma (C'')-(p,q)$ and $C'''=\sigma (C''')-(p,q)$ if and only if
$(i-p,j-q)$ is a southwest angle of $\Gamma$  southwest of $(i,j)$
such that $-i=i-p$ and $j-q-1\neq -j$.
\Item[d)] Otherwise $C\setminus \{(-i,-j),(i-p,j-q)\}=C'$ is a connected 
component of $\Delta_{p,q}'$ such that $C'=\sigma (C')-(p,q)$.
\end{Itemize}
\end{Lemma}

All three lemmas are easy consequences of Lemma 3.3, of the central symmetry
of $\Gamma$ and the fact that $(i,j)$ is chosen to be the northeast
corner with $i$ minimal.
\vskip0.5em

It follows from these three lemmas that $e_{p,q}(\Gamma ,\Gamma )$ and
$e_{p,q}(\Gamma' ,\Gamma')$ are equal unless $(p,q)$ satisfies one of the
following conditions:

\begin{Itemize}
\Item[1.] $(i-p,j-q)\neq (-i,-j)$ is a southwest corner of $\Gamma$ 
southwest of $(i,j)$, in which case 
$e_{p,q}(\Gamma ,\Gamma )=e_{p,q}(\Gamma',\Gamma')+1$
by Lemma 3.7 a). 

\Item[2.] $(i-p,j-q)$ is a southwest angle of $\Gamma$ southwest
of $(i,j)$ such that $i-p\neq -i$, in which case 
$e_{p,q}(\Gamma ,\Gamma)=e_{p,q}(\Gamma',\Gamma')-1$ by Lemma 3.7 c).

\Item[3.] $(i-p,j-q)$ is a southwest angle of $\Gamma$ southwest of $(i,j)$
such that $i-p=-i$ and $j-q\neq -j+1$, in which case
$e_{p,q}(\Gamma ,\Gamma)=e_{p,q}(\Gamma',\Gamma') -1$ by Lemma 3.8 b).

\Item[4.] $(i-p,j-q)=(-i,-j)$ with $i\geq 0$, in which case 
\ItemItem[i)] $e_{p,q}(\Gamma ,\Gamma)=e_{p,q}(\Gamma' ,\Gamma')$ if
$2(i+j)$ is even by Lemma 3.8 a);
\ItemItem[ii)] $e_{p,q}(\Gamma ,\Gamma)=e_{p,q}(\Gamma' ,\Gamma')$+1 if
$2(i+j)$ is odd by Lemma 3.8 a).

\Item[5.] $(i-p,j-q)=(-i+1,-j+1)$ with $i>0$, in which case
\ItemItem[i)] $e_{p,q}(\Gamma ,\Gamma)=e_{p,q}(\Gamma' ,\Gamma')-1$ if
$2(i+j)$ is even by Lemma 3.8 b);
\ItemItem[ii)] $e_{p,q}(\Gamma ,\Gamma)=e_{p,q}(\Gamma' ,\Gamma')$ if
$2(i+j)$ is odd by Lemma 3.8 b).

\Item[6.] $(i-p,j-q)=(-i+1,-j)$ with $i>0$, in which case
\ItemItem[i)] $e_{p,q}(\Gamma ,\Gamma)=e_{p,q}(\Gamma' ,\Gamma')+1$ if
$2(i+j)$ is even by Lemma 3.8 a);
\ItemItem[ii)] $e_{p,q}(\Gamma ,\Gamma)=e_{p,q}(\Gamma' ,\Gamma')$ if
$2(i+j)$ is odd by Lemma 3.8 a).

\Item[7.] $(i-p,j-q)=(-i,-j+1)$ is a southwest angle of $\Gamma$ southwest
of $(i,j)$, in which case
\ItemItem[i)] $e_{p,q}(\Gamma ,\Gamma)=e_{p,q}(\Gamma' ,\Gamma')$ if
$2(i+j)$ is even by Lemma 3.8 b);
\ItemItem[ii)] $e_{p,q}(\Gamma ,\Gamma)=e_{p,q}(\Gamma' ,\Gamma')-1$ if
$2(i+j)$ is odd by Lemma 3.8 b).

\Item[8.] $(i-p,j-q)=(-i,-j+1)$ is not a southwest angle of $\Gamma$ southwest
of $(i,j)$ and $i\geq 0$, in which case
\ItemItem[i)] $e_{p,q}(\Gamma ,\Gamma)=e_{p,q}(\Gamma' ,\Gamma')+1$ if
$2(i+j)$ is even by Lemma 3.8 a);
\ItemItem[ii)] $e_{p,q}(\Gamma ,\Gamma)=e_{p,q}(\Gamma' ,\Gamma')$ if
$2(i+j)$ is odd by Lemma 3.8 a).
\end{Itemize}
\vskip0.5em

Now if $i<0$, then only 1 and 2 apply, thus we deduce as in Theorem 3.2
that
$$
\sum_{p,q\geq 0}e_{p,q}(\Gamma ,\Gamma )=
\sum_{p,q\geq 0}e_{p,q}(\Gamma',\Gamma')+N_{c}-N_{a}
$$
where $N_{c}$ (resp. $N_{a}$) is the number of southwest corners (resp. angles)
southwest of $(i,j)$.

Now if $i=0$, then 5 and 6 do not apply. We have therefore
$$
\sum_{p,q\geq 0}e_{p,q}(\Gamma ,\Gamma )=
\sum_{p>0,q\geq 0}e_{p,q}(\Gamma ,\Gamma )+
\sum_{q\geq 0}e_{0,q}(\Gamma ,\Gamma ).
$$
There are 4 cases: $2(i+j)$ is even or odd; $(-i,-j+1)$ is or is not
a southwest angle of $\Gamma$. By applying 3,4,7 and 8 in the appropriate
case, we have that:
$$
\sum_{q\geq 0}e_{0,q}(\Gamma ,\Gamma )=
\sum_{q\geq 0}e_{0,q}(\Gamma' ,\Gamma' )
$$
since $\Gamma$ is not rectangular, so there is always a southwest angle
above $(-i,-j)$. It follows that we have again
$$
\sum_{p,q\geq 0}e_{p,q}(\Gamma ,\Gamma )=
\sum_{p,q\geq 0}e_{p,q}(\Gamma',\Gamma')+N_{c}-N_{a}.
$$
Finally if $i>0$, then 5 and 6 apply also. Again there are the same 4 cases to 
consider, and we deduce that
$$
\sum_{0\leq p\leq 2i,q\geq 0}e_{p,q}(\Gamma ,\Gamma )=
\sum_{0\leq p\leq 2i,q\geq 0}e_{p,q}(\Gamma',\Gamma')
$$
since for these $(p,q)$, $e_{p,q}(\Gamma ,\Gamma )=e_{p,q}(\Gamma',\Gamma')$
unless 3,4,5,6,7 or 8 is satisfied (note that $(-i,-j)$ is the first
southwest corner of $\Gamma$ southwest of $(i,j)$ because $i>0$ and $(i,j)$
is chosen to be the leftmost northeast corner of $\Gamma$).
Hence, we have again
$$
\sum_{p,q\geq 0}e_{p,q}(\Gamma ,\Gamma )=
\sum_{p,q\geq 0}e_{p,q}(\Gamma',\Gamma')+N_{c}-N_{a}.
$$ 
By induction, we have our result.

We are therefore left with the case where $\Gamma'$ is not connected.
In this case, let $(i,j)$ be the leftmost southwest corner of $\Gamma$.
Then $\Gamma$ is of the following form:

%
% Basic example 5
%
$$
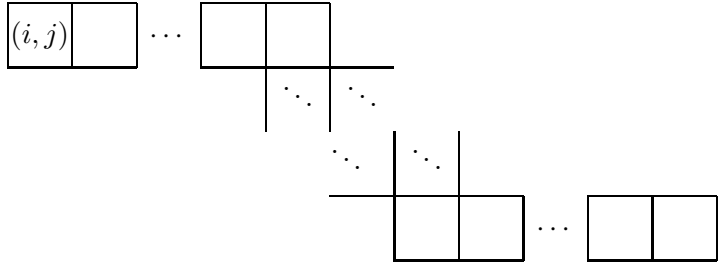
\vbox{\offinterlineskip
\halign{#&#&#&#&#&#&#&#&#&#&#&#&#\cr
\omit\hrulefill&\omit\hrulefill&&\omit\hrulefill&\omit\hrulefill&\cr
\boite{$(i,j)$}&\boite{}&\boite{$\cdots$}&\boite{}&\boite{}&\boite{}\cr
\omit\hrulefill&\omit\hrulefill&&\omit\hrulefill&\omit\hrulefill&\omit\hrulefill\cr
&&&&\boite{$\ddots$}&\boite{$\ddots$}\cr
&&&&&$\ddots$&\boite{$\ddots$}&\boite{}\cr
&&&&&\omit\hrulefill&\omit\hrulefill&\omit\hrulefill&&\omit\hrulefill&\omit\hrulefill\cr
&&&&&&\boite{}&\boite{}&\boite{$\cdots$}&\boite{}&\boite{}&\boite{}\cr
&&&&&&\omit\hrulefill&\omit\hrulefill&&\omit\hrulefill&\omit\hrulefill\cr
}
}
$$

A similar argument as in Theorem 3.2 reduces $\Gamma$ to a horizontal
diagram which is a simple verification.
\end{proof}

\begin{Theorem}
Let $\Gamma,\Gamma'$ be two centrally symmetric skew diagrams
such that $\Gamma\cap \Gamma'=\{ (0,0)\}\neq \Gamma'$. Then  
$\Card \bigcup_{p,q\geq 0, (p,q)\neq (0,0)} 
\E_{p,q}(\Gamma ,\Gamma' ) \geq 1$.
Equality holds if and only if $\Gamma$ and $\Gamma'$ satisfy one of the
following two conditions:
\begin{Itemize}
\Item[a)] $\Gamma =\{ (0,0)\}$ and $\Gamma'$ has exactly one southwest corner
$(i,j)$ such that $i,j\leq 0$;  
\Item[b)] $\Gamma$ and $\Gamma'$ are both rectangular, {\it i.e.} they both
have exactly one southwest corner and one northeast corner.
\end{Itemize}
\end{Theorem}
\begin{proof}
Let us first suppose that $\Gamma =\{ (0,0)\}$. Then $\Gamma'$ has at least
a northeast corner $(p,q)$ such that $p,q\geq 0$, $(p,q)\neq (0,0)$ 
since $\Gamma'$ is a skew diagram of cardinal $>1$ containing $(0,0)$. Thus
the pair $\{(0,0)\}$, $\{(-p,-q)\}$ is in $\E_{p,q}(\Gamma ,\Gamma' )$. 

We can therefore suppose that $\Gamma$ and $\Gamma'$ are skew diagrams of 
cardinal $> 2$ and $\Gamma\cap \Gamma'=\{(0,0)\}$. 
There exist $p,q> 0$ such that, exchanging the roles of $\Gamma$ and 
$\Gamma'$ if necessary, 
$(-p,0)$ be a southwest corner of $\Gamma$ and
$(0,q)$ be a northeast corner of $\Gamma'$. By central
symmetry, $(0,-q)\in\Gamma'$ and the pair $\{ (-p,0)\} , \{(0,-q)\}$
is in $\E_{p,q}(\Gamma,\Gamma')$. So the first part of the
theorem is proved.
\vskip1em

If $\Gamma =\{ (0,0)\}$, then it is clear that equality holds if and only
if condition a) is satisfied. So we can assume that the cardinality
of $\Gamma$ and $\Gamma'$ are $>2$.

If they are both rectangular, then one is horizontal and the other vertical.
It is clear that equality holds. Now let $\Gamma$ be non rectangular and
that $\Gamma$ contains $(-1,0)$ (the case $\Gamma$ contains $(0,-1)$ is 
similar). Then $\Gamma'$ contains $(0,-1)$.

Let $(0,q)$ be a northeast
corner of $\Gamma'$, then $q>0$ and $(0,q-1)\in\Gamma'$.
Let $(-p,0)$ be a southwest corner of $\Gamma$ with $p>0$.
Then either $(-p-1,1)$ is not in $\Gamma$ or $\Gamma$ has
a southwest corner $(-r,1)$ with $r>0$. 

In the first case, the pair $\{(-p,0),(-p,1)\}$, $\{(0,-q),(0,-q+1)\}$ is
in $\E_{p,q}(\Gamma,\Gamma')$. In the second case, the pair
$(-r,1),(0,-q)$ is in $\E_{r,q-1}(\Gamma,\Gamma')$. Since the pair
$(-p,0),(0,-q)$ is in $\E_{p,q}(\Gamma,\Gamma')$ in both cases,
we are done.
\end{proof}

%
% Counterexample
% 
\begin{Example}
Let $\Gamma = \{(-1,0),(0,0),(1,0)\}$ and
$\Gamma' = \{(-1,1),(0,1)$, $(0,0),$ $(0,-1),$ $(1,-1)\}$ be centrally
symmetric skew diagrams. 
$$
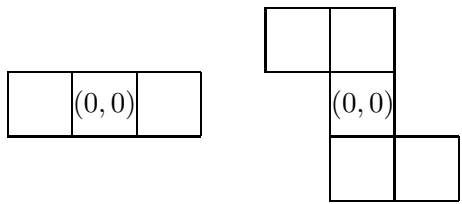
\vbox{\offinterlineskip
\halign{#&#&#&#&#&#&#&#&#&#&#&#&#&#&#&#\cr
&&&&&&\omit\hrulefill&\omit\hrulefill&&\cr
&&&&&&\boite{}&\boite{}&\boite{}&\cr
&\omit\hrulefill&\omit\hrulefill&\omit\hrulefill&&&\omit\hrulefill&\omit\hrulefill&&\cr
&\boite{}&\boite{$(0,0)$}&\boite{}&\boite{}&&&\boite{$(0,0)$}&\boite{}&\cr
&\omit\hrulefill&\omit\hrulefill&\omit\hrulefill&&&&\omit\hrulefill&\omit\hrulefill&\cr
&&&&&&&\boite{}&\boite{}&\boite{}\cr
&&&&&&&\omit\hrulefill&\omit\hrulefill&\cr
}
}
$$
Then
$\E_{1,1}(\Gamma ,\Gamma')$ contains the unique pair $(-1,0),(0,-1)$
and $\E_{0,1}(\Gamma ,\Gamma')$ contains the unique pair
$\{ (-1,0),(0,0)\},\{ (0,-1),(1,-1)\}$. One verifies easily that
the cardinality of the set $\bigcup_{p,q\geq 0}\E_{p,q}(\Gamma,\Gamma')$ is
$2$. 
\end{Example}

\begin{Proposition}
Let $\Gamma$, $\Gamma'$ be centrally symmetric skew diagrams of different
types. Then $\bigcup_{p,q\geq 0, (p,q)\neq (0,0)}
\E_{p,q}(\Gamma ,\Gamma')$ is non empty.
\end{Proposition}
\begin{proof}
A centrally symmetric skew diagram has always a southwest corner $(i,j)$
such that $i,j\leq 0$. Let $(i,j)$ (resp. $(k,l)$) be such a southwest
corner in $\Gamma$ (resp. $\Gamma'$). Then it is clear
that the pair $(i,j), (k,l)$ is in $\E_{-i-k,-j-l}(\Gamma ,\Gamma')$ and
we are done.
\end{proof}

\section{Conditions (Y), (R) and Near rectangular skew diagrams}

Let $\Gamma$ be a skew diagram (resp. centrally symmetric skew diagram)
and set 

\begin{Itemize}

\Item[$\bullet$] $\E (\Gamma ):=\bigcup_{(p,q)\in\Z^{2}}\E_{p,q}(\Gamma )$
(resp. $\E (\Gamma ,\Gamma ):=\bigcup_{(p,q)\in\Z^{2}}\E_{p,q}(\Gamma ,\Gamma)$)

\Item[$\bullet$] and 
$\E_{+}(\Gamma ):=\bigcup_{p,q\geq 0}\E_{p,q}(\Gamma )$
(resp.
$\E_{+} (\Gamma,\Gamma ):=\bigcup_{p,q\geq 0}\E_{p,q}(\Gamma ,\Gamma)$). 
\end{Itemize}

From the previous section, we have that 
$\Card \E_{+}(\Gamma )= \Card \Gamma -1$
(resp. $\Card \E_{+}(\Gamma,\Gamma )= [\Card \Gamma /2]$).
In this section, we shall study skew diagrams (resp. centrally
symmetric skew diagrams) 
such that the cardinality of 
$\E (\Gamma )$ (resp. $\E (\Gamma,\Gamma )$) 
is less than or equal to $1+\Card \E_{+}(\Gamma )$.
(resp. $1+ \Card \E_{+}(\Gamma,\Gamma )$).

\begin{Definition}
We say that a skew diagram $\Gamma$ satisfies (Y) if it has either
exactly $1$ southwest corner or exactly $1$ northeast corner.
\end{Definition}

So in the terminology of \cite{[G]}, $\Gamma$ satisfies (Y) if and only
if $\Gamma$ is a Young diagram or a minus Young diagram, while in the
terminology of \cite{[EP]}, this is equivalent to that $\Gamma$ is either
sw-Young or ne-Young.
 
\begin{Proposition}
Let $\Gamma$ be a skew diagram. Then 

\begin{Itemize}
\Item[a)] $\E (\Gamma )=\E_{+} (\Gamma )$ 
if and only if $\Gamma$ satisfies (Y).
\Item[b)] if $\Gamma$ does not satisfy (Y), then
$\Card \E (\Gamma )\geq 2+\Card \E_{+} (\Gamma )$.
\end{Itemize}
\end{Proposition}
\begin{proof}
First note that $\Gamma$ has exactly one southwest corner is equivalent to
$\sigma(\Gamma )$ has exactly one northeast corner, and that if
$C\in \E_{p,q}(\Gamma )$, then $\sigma (C)-(p,q)\in \E_{p,q}(\sigma (\Gamma ))$.
So we can suppose that $\Gamma$ has exactly one southwest corner $(i,j)$.
Now any skew subdiagram of $\Gamma$ contains $(i,j)$. So $\E_{p,q}(\Gamma )$
is empty if $p<0$ or $q<0$, and the ``if'' part follows.

Now suppose that $\Gamma$ does not satisfy (Y), and therefore there 
are at least $2$ southwest corners and at least $2$ northeast corners.
Let the bottom row of $\Gamma$ be $(i,j),\cdots ,(i+r,j)$
and the top row of $\Gamma$ be $(k,l),\cdots ,(k+s,l)$ where $r,s\geq 0$.
Then $\Gamma$ does not satisfy (Y) implies that $k<i$, $j<l$
and $i+r>k+s$.

Now let $t=\min (r,s)$, then $k+s-t<i$. It follows that
the skew diagram $\{ (i,j),\cdots ,(i+t,j)\}$ is in 
$\E_{k+s-t-i,l-j}(\Gamma )$. So part a) follows
since $k+s-t-i < 0$.

Finally, if $\Gamma$ does not satisfy (Y), then the above can also
be applied to the leftmost column and the rightmost column. Thus b)
follows.
\end{proof}

\begin{Definition}
A centrally symmetric skew diagram $\Gamma$ is called {\bf near rectangular}
if it has exactly $2$ southwest corners $(i,j)$, $(k,l)$ with $i<k$ such that 
one of the following conditions is satisfied:
\begin{Itemize}
\Item[a)] $\Gamma \setminus \{ (i,j),(-i,-j)\}$ is rectangular;
\Item[b)] $\Gamma \setminus \{ (k,l),(-k,-l)\}$ is rectangular;
\Item[c)] $(i,j)=(k-1,l+1)$.
\end{Itemize}
Their standard forms are illustrated as follows:

%
% NR
%
$$
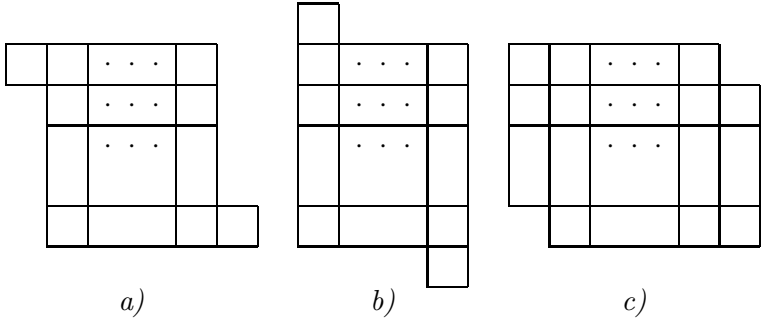
\vbox{\offinterlineskip
\halign{#&#&#&#&#&#&#&#&#&#&#&#&#&#&#&#&#&#&#&#&#&#&#&#&#\cr
&&&&&&&\omit\hrulefill\cr
&&&&&&&\sboite{}&\sboite{}\cr
\omit\hrulefill&\omit\hrulefill&\omit\hrulefill&\omit\hrulefill&&&
&\omit\hrulefill&\omit\hrulefill&\omit\hrulefill&&&
\omit\hrulefill&\omit\hrulefill&\omit\hrulefill&\omit\hrulefill&\cr
\sboite{}&\sboite{}&\sboite{$\cdot$}$\cdot$\hboite{$\cdot$}&\sboite{}&\sboite{}&&
&\sboite{}&\sboite{$\cdot$}$\cdot$\hboite{$\cdot$}&\sboite{}&\sboite{}&&
\sboite{}&\sboite{}&\sboite{$\cdot$}$\cdot$\hboite{$\cdot$}&\sboite{}&\sboite{}\cr
\omit\hrulefill&\omit\hrulefill&\omit\hrulefill&\omit\hrulefill&&&
&\omit\hrulefill&\omit\hrulefill&\omit\hrulefill&&&
\omit\hrulefill&\omit\hrulefill&\omit\hrulefill&\omit\hrulefill&\omit\hrulefill&\cr
&\sboite{}&\sboite{$\cdot$}$\cdot$\hboite{$\cdot$}&\sboite{}&\sboite{}&&
&\sboite{}&\sboite{$\cdot$}$\cdot$\hboite{$\cdot$}&\sboite{}&\sboite{}&&
\sboite{}&\sboite{}&\sboite{$\cdot$}$\cdot$\hboite{$\cdot$}&\sboite{}&\sboite{}&\sboite{}\cr
&\omit\hrulefill&\omit\hrulefill&\omit\hrulefill&&&
&\omit\hrulefill&\omit\hrulefill&\omit\hrulefill&&&
\omit\hrulefill&\omit\hrulefill&\omit\hrulefill&\omit\hrulefill&\omit\hrulefill&\cr
&\sboite{}&\sboite{$\cdot$}$\cdot$\hboite{$\cdot$}&\sboite{}&\sboite{}&&
&\sboite{}&\sboite{$\cdot$}$\cdot$\hboite{$\cdot$}&\sboite{}&\sboite{}&&
\sboite{}&\sboite{}&\sboite{$\cdot$}$\cdot$\hboite{$\cdot$}&\sboite{}&\sboite{}&\sboite{}\cr
&\sboite{}&\sboite{}&\sboite{}&\sboite{}&&
&\sboite{}&\sboite{}&\sboite{}&\sboite{}&&
\sboite{}&\sboite{}&\sboite{}&\sboite{}&\sboite{}&\sboite{}\cr
&\omit\hrulefill&\omit\hrulefill&\omit\hrulefill&\hrulefill&&
&\omit\hrulefill&\omit\hrulefill&\omit\hrulefill&&&
\omit\hrulefill&\omit\hrulefill&\omit\hrulefill&\omit\hrulefill&\omit\hrulefill&\cr
&\sboite{}&\sboite{}&\sboite{}&\sboite{}&\sboite{}&
&\sboite{}&\sboite{}&\sboite{}&\sboite{}&&
&\sboite{}&\sboite{}&\sboite{}&\sboite{}&\sboite{}\cr
&\omit\hrulefill&\omit\hrulefill&\omit\hrulefill&\hrulefill&&
&\omit\hrulefill&\omit\hrulefill&\omit\hrulefill&&&
&\omit\hrulefill&\omit\hrulefill&\omit\hrulefill&\omit\hrulefill&\cr
&&&&&&&&&\sboite{}&\sboite{}\cr
&&&&&&&&&\omit\hrulefill\cr
&\cr
&&\hfil a)\hfil &&&&&&\hfil b)\hfil&&&&&&\hfil c)\hfil\cr
}
}
$$

We say that a skew diagram $\Gamma$ satisfies (R) if
either $\Gamma$ is rectangular or $\Gamma$ is non integral and 
near rectangular.
\end{Definition}

\begin{Proposition}
Let $\Gamma$ be a centrally symmetric skew diagram. Then
$\E (\Gamma ,\Gamma )=\E_{+}(\Gamma ,\Gamma )$
if and only if $\Gamma$ satisfies (R).
\end{Proposition}
\begin{proof}
Let us suppose that $\Gamma$ satisfies (R). If $\Gamma$ is rectangular,
then we have our result by Proposition 4.2. So let us suppose that
$\Gamma$ is near rectangular with southwest corners $(i,j)$, $(k,l)$
with $i<k$. 

If $\Gamma$ satisfies condition a) of Definition 4.3, 
then $k=i+1$ and $l=-j$. We observe that we have only to find
a necessary and sufficient condition for 
$\{ (i,j)\}$ to be in $\E_{-2i,-2j}(\Gamma ,\Gamma )$
or $\{ (i+1,-j),\cdots ,(-i,-j)\}$ to be in $\E_{-1,2j}(\Gamma ,\Gamma )$. 
This is so if and only if $2(i+j)$ or $2j-1$ is odd.
But this is equivalent to saying that $\Gamma$ is not non integral.

If $\Gamma$ satisfies condition b) of Definition 4.3,
then a similar argument gives our result.

Now if $\Gamma$ satisfies condition c) of Definition 4.3,
then $(k,l)=(i+1,j-1)$. We now have to find a necessary and
sufficient condition for
$\{ (i,j),\cdots (i,-j+1)\}$ to be in $\E_{-2i,-1}(\Gamma ,\Gamma )$
or $\{ (i+1,j-1),\cdots ,(-i,j-1)\}$ to be in 
$\E_{-1,-2j+2}(\Gamma ,\Gamma )$.
This is so if and only if $2i+1$ or $2j+1$ is odd.
Again this is equivalent to saying that $\Gamma$ is not non integral.

We have therefore proved the ``if'' part.

To prove the ``only if'' part, it suffices therefore
to prove that if $\E (\Gamma ,\Gamma )
=\E_{+}(\Gamma ,\Gamma )$, then
$\Gamma$ is either rectangular or near rectangular.

Let $(i,j),\cdots ,(i+r,j)$ be the top row of $\Gamma$ with $-i,j,r\geq 0$,
then $(-i-r,-j),\cdots ,(-i,-j)$ is the bottom row of $\Gamma$.
Note that $i\leq -i-r$ since $\Gamma$ is a skew diagram. 
If $r \geq 1$, then either 
$\{ (-i-r,-j),\cdots ,(-i,-j)\}$ is in $\E_{2i+r,2j}(\Gamma ,\Gamma )$
or $\{ (-i-r,j),\cdots ,(-i-1,j)\}$ is in $\E_{2i+r+1,2j}(\Gamma ,\Gamma )$
since either $2i+r+2j$ or $2i+r+2j+1$ is odd. 
So it follows that if $r \geq 1$, then $\E_{p,q}(\Gamma ,\Gamma )$ is
non empty for some $p,q$ with $p<0$, $q>0$ if $2i+r+1 < 0$, 
or equivalently, $i+1 < -i-r$.  

Now let $(k,l),\cdots ,(k,l+s)=(i,j)$ be the leftmost column
of $\Gamma$ with $-k,-l,s\geq 0$, then
$(-k,-l-s),\cdots ,(-k,-l)$ is the rightmost column of $\Gamma$.
Note that $l\geq -l-s$ since $\Gamma$ is a skew diagram.
A similar argument shows that if $s\geq 1$, then 
$\E_{p,q}(\Gamma ,\Gamma )$ is non empty for some $p,q$ with 
$q<0$, $p>0$ if $-l-s<l-1$.

We are therefore left with four possibilities: 
$r=0$ and $s=0$; $r=0$ and $s> 0$; $r>0$ and $s=0$;
$r>0$ and $s>0$.

In the first case, we have that $\Gamma =\{ (0,0)\}$ which is rectangular
and the result follows.

In the second case, we have that $l-1\leq -l-s$.
Since $l\geq -l-s$, we have either $l=-l-s$ or $l-1=-l-s$. If $l=-l-s$, 
then $\Gamma$ is just a vertical diagram which is rectangular. 
If $l-1=-l-s$, then $\Gamma$ is near rectangular satisfying condition b)
of Definition 4.3.

The third case is similar to the second yielding condition a) of
Definition 4.3. The last case implies
that either $l=-l-s$, $i=-i-r$ or $l-1=-l-s$, $i+1=-i-r$.
In the former, $\Gamma$ is rectangular. In the latter,
$\Gamma$ is near rectangular satisfying condition c) of Definition 4.3.

Our proposition now follows.
\end{proof}

\begin{Proposition}
Let $\Gamma$ be a centrally symmetric skew diagram. Then
$\Card \E (\Gamma,\Gamma)=1+\Card \E_{+}(\Gamma,\Gamma )$ if 
and only if one of the following conditions is satisfied:
\begin{Itemize}
\Item[i)] $\Gamma$ is semi-integral and near rectangular;
\Item[ii)] $\Gamma$ is integral and near rectangular of type a) or b).
\end{Itemize}
\end{Proposition}
\begin{proof}
The ``if'' part is just a direct verification using the computations
in the proof of Proposition 4.4. 
Let us show why type c) is not allowed for $\Gamma$ integral and 
leave the other verifications to the reader. 
Let $\Gamma$ be integral and near rectangular
and $(i,j)$, $(i+1,j-1)$ be the two distinct southwest corners with $i<k$.
Then $2i+1$ and $2j+1$ are both odd, so by the computations
in the proof of Proposition 4.4, there are 2 elements of
$\E (\Gamma,\Gamma )$ which are not in $\E_{+}(\Gamma ,\Gamma )$.
Note that if $\Gamma$ is near rectangular and non integral, then
$\E (\Gamma,\Gamma )=\E_{+}(\Gamma,\Gamma )$. 

To prove the ``only if'' part, it suffices to prove that if 
$\Card \E (\Gamma,\Gamma)=1+\Card \E_{+}(\Gamma,\Gamma )$, then
$\Gamma$ is near rectangular.
Let us suppose that $\Gamma$ is not near rectangular. 

Let $(i,j),\cdots ,(i+r,j)$ be the top row of $\Gamma$ with $-i,j,r\geq 0$,
then $(-i-r,-j),\cdots ,(-i,-j)$ is the bottom row of $\Gamma$.
Note that $i\leq -i-r$ since $\Gamma$ is a skew diagram. 

Let $(k,l),\cdots ,(k,l+s)=(i,j)$ be the leftmost column
of $\Gamma$ with $-k,-l$, $s\geq 0$, then
$(-k,-l-s),\cdots ,(-k,-l)$ is the rightmost column of $\Gamma$.
Note that $l\geq -l-s$ since $\Gamma$ is a skew diagram.

Recall from the proof of Proposition 4.4 that if $r \geq 1$ and $i+1<-i-r$,
then $\E_{p,q}(\Gamma ,\Gamma )$ is non empty for some $p,q$ with $p<0$, $q>0$;
and if $s\geq 1$ and $-l-s<l-1$, then 
$\E_{p,q}(\Gamma ,\Gamma )$ is non empty for some $p,q$ with 
$p>0$, $q<0$.

As in Proposition 4.4, we have four cases to consider: 
$r=s=0$; $r=0$, $s>0$; $s=0$, $r>0$; $r,s> 0$.

In the first case, $\Gamma$ is rectangular.
So $\Card \E (\Gamma,\Gamma)=\Card \E_{+}(\Gamma,\Gamma )$.

In the second case, if $s=-2l$, then $\Gamma$ is rectangular
and $\Card \E (\Gamma,\Gamma)=\Card \E_{+}(\Gamma,\Gamma )$. 
Since $\Gamma$ is not near rectangular, $-l-s < l-1$ and $k<0$. 
So $\E_{p,q}(\Gamma ,\Gamma )$ is non empty for some $p,q$ with 
$p>0$ and $q<0$ if $-l-s<l-1$.

Let $(u,v)$ be a southwest corner distinct from
$(k,l)$ and $(-i,-j)$. Then clearly, $i=k< u< -i$.
It follows that $\{ (u,v)\}, \{(-i,-j)\}$ is a pair
in $\E_{i-u,j-v}(\Gamma,\Gamma )$. Therefore 
$\Card \E (\Gamma,\Gamma )\geq \Card \E_{+}(\Gamma ,\Gamma )+2$.

So let us assume that $(k,l)$ and $(-i,-j)$ are
the only southwest corners of $\Gamma$. Then since $-l-s < l-1$, 
$(k+1,l+s-1)\not \in\Gamma$ and by central symmetry
$(-i-1,-j+1)\not \in\Gamma$; otherwise, $(-i-1,-j+1)$ is a
southwest corner distinct from $(k,l)$ and $(-i,-j)$.
It follows that either $\{(-i,-j)\}$ is in 
$\E_{k+i,l+s+j}(\Gamma,\Gamma)$ or $\{(-i,-j),(-i,-j+1)\}$ is in 
$\E_{k+i,l+s+j-1}(\Gamma,\Gamma)$, where $k+i < 0$. Thus
$\Card \E (\Gamma,\Gamma )\geq \Card \E_{+}(\Gamma ,\Gamma )+2$.

The third case is similar to the second. Finally, in the
last case, we have either $i+1\geq -i-r$, $-l-s<l-1$ or
$i+1<-i-r$, $-l-s\geq l-1$. Without loss of generality,
we can assume that $i+1\geq -i-r$, $-l-s<l-1$. Then
$i+r=-i$ or $i+r=-i-1$. Since $-l-s<l-1$, $\Gamma$ is not 
rectangular and therefore we can not have $i+r=-i$. 
So $i+r=-i-1$ and $\Gamma$ has exactly two southwest corners
$(k,l)=(i,j-s)$ and $(k+1,-l-s)=(i+1,-j)$.

Now $\E_{p,q}(\Gamma ,\Gamma )$ is non empty for some $p,q$ with 
$p>0$ and $q<0$ since $-l-s<l-1$. 
Further, $-l-s<l-1$ implies that $(i,-j+1)\not \in\Gamma$, therefore
the rectangular diagram formed by the bottom two rows
$\{ (i+1,-j),\cdots ,(-i,-j),(i+1,-j+1), \cdots ,(-i,-j+1)\}$ is
a skew subdiagram of $\Gamma$. It follows that
either the bottom row $\{ (i+1,-j),\cdots ,(-i,-j)\}$ is
in $\E_{k-i-1,l+s+j}(\Gamma,\Gamma )$ or 
$\{ (i+1,-j),\cdots ,(-i,-j),(i+1,-j+1), \cdots ,(-i,-j+1)\}$
is in $\E_{k-i-1,l+s-1+j}(\Gamma ,\Gamma )$ where $k-i-1=-1<0$.
Again we conclude that 
$\Card \E (\Gamma,\Gamma )\geq \Card \E_{+}(\Gamma ,\Gamma )+2$.

We have therefore proved that if $\Gamma$ is not 
near rectangular, then we have
$\Card \E (\Gamma,\Gamma )\neq \Card \E_{+}(\Gamma ,\Gamma )+1$. 
Our proposition now follows.
\end{proof}

\section{Panyushev's conjecture}

We shall now apply the results of the previous sections to 
answer Panyushev's conjecture which states that distinguished 
nilpotent pairs are wonderful. 

Let us denote by $Z_{+}(\ee )$ the direct sum 
$\bigoplus_{p,q\geq 0} Z_{p,q}(\ee )$.

\begin{Theorem}
The dimension of the positive quadrant $Z_{+}(\ee )$ is greater than or equal 
to the rank of $\ggoth$. Equality holds if and only if one of the following
conditions is satisfied: 
\begin{Itemize}
\Item[a)] $\ggoth$ is of type A;
\Item[b)] $\ggoth$ is of type B and the
associated pair of centrally symmetric skew diagrams (see Theorem 1.5)
$(\Gamma^{1},\Gamma^{2})$ satisfies $\Gamma^{2}=\emptyset$.
\Item[c)] $\ggoth$ is of type C and the
associated pair of centrally symmetric skew diagrams (see Theorem 1.5)
$(\Gamma^{1},\Gamma^{2})$ satisfies
$\Gamma^{1}=\emptyset$ or $\Gamma^{2}=\emptyset$.
\Item[d)] $\ggoth$ is of type $D$ and the associated triple of
centrally symmetric skew diagrams (see Theorem 1.5) 
$(\Gamma^{1},\Gamma^{2},\Gamma^{3})$ is such that either 
$\Gamma^{2}=\Gamma^{3}=\emptyset$ or $\Gamma^{1}$ is empty
and the pair $(\Gamma^{2},\Gamma^{3})$ satisfies the conditions a) or b) of 
Theorem 3.9.
\end{Itemize}
\end{Theorem}
\begin{proof}
This is a direct consequence of the Theorem 2.9, Theorems 3.2,
3.5, 3.9 and Proposition 3.11.
\end{proof}
 
\begin{Definition} {\rm \cite{[P1]}}
Recall that a nilpotent pair $\ee =(e_{1},e_{2})$ is called {\bf wonderful} 
if $\dim \bigoplus_{p,q\in\Z, p,q\geq 0}Z_{p,q}(\ee )=\rank \ggoth$.
\end{Definition}

It follows from Theorems 2.9, 3.2, 3.5 and 3.9 that:

\begin{Corollary}
A distinguished nilpotent pair $\ee$ in a classical simple Lie algebra 
$\ggoth$ is wonderful if and only if one of the following
conditions is satisfied:
\begin{Itemize}
\Item[a)] $\ggoth$ is of type A,B or C;
\Item[b)] $\ggoth$ is of type $D$ and the associated triples of
centrally symmetric skew diagrams (see Theorem 1.5) 
$(\Gamma^{1},\Gamma^{2},\Gamma^{3})$ is such that the pair 
$(\Gamma^{2},\Gamma^{3})$ satisfies the conditions a) or b) of 
Theorem 3.9.
\end{Itemize}
\end{Corollary}
\begin{proof}
It suffices to observe that Proposition 3.11 does not apply here
since we are only interested in $p,q\in\Z$. 
\end{proof}

So distinguished nilpotent pairs are always wonderful in types A,
B and C, but they are not so in type D. 
Therefore Panyushev's
conjecture is not true in type D.

\begin{Example}
An example of a distinguished nilpotent pair which is not
wonderful is the one described in Example 3.10 which corresponds
to a distinguished nilpotent pair in $D_{4}$. One verifies easily that
the positive quadrant has dimension $5$.
\end{Example}

\section{Classification of principal and
almost principal nilpotent pairs}

By applying the Theorem 5.1 and Propositions 4.2, 4.4, we recover easily
the classification of principal nilpotent pairs.

\begin{Corollary}
Let $\ee$ be a distinguished nilpotent pair in a classical simple
Lie algebra $\ggoth$. Then  
$\dim Z(\ee )=\rank \ggoth$ if and only if
\begin{Itemize}
\Item[a)] $\ggoth$ is of type A and the associated skew diagram $\Gamma$
satisfies (Y);
\Item[b)] $\ggoth$ is of type B and the associated pair of 
centrally symmetric skew diagrams $(\Gamma^{1},\emptyset )$ satisfies 
(R);
\Item[c)] $\ggoth$ is of type C and the associated pair of
centrally symmetric skew diagrams $(\Gamma^{1},\emptyset )$ (resp.
$(\emptyset ,\Gamma^{2})$) satisfies (R);
\Item[d)] $\ggoth$ is of type $D$ and the associated triples of
centrally symmetric skew diagrams
$(\emptyset ,\Gamma^{2},\Gamma^{3})$ (resp. $(\Gamma^{1},\emptyset,\emptyset )$)
satisfies (R).
\end{Itemize}
\end{Corollary}

The above corollary is exactly the classification of principal nilpotent 
pairs in classical simple Lie algebra given in \cite{[EP]}. 

\begin{Remarks}
Note that there exist non principal distinguished nilpotent pairs satisfying 
$Z (\ee )=Z_{+}(\ee )$. By the results of Sections 3 and 4, such a nilpotent
pair must be of type other than A and the corresponding skew diagrams
must all satisfy (R).

For example, in type 
B and C, take any pair $(\Gamma^{1},\Gamma^{2})$ such that they are both 
non empty and rectangular. Then one can verify easily that
$Z (\ee )=Z_{+}(\ee )$. 

Also in type B, take any pair $(\Gamma^{1},\Gamma^{2})$ such that $\Gamma^{1}$
is rectangular with southwest corner $(i,j)$ and $\Gamma^{2}$ is near 
rectangular satisfying for all $(k,l)\in \Gamma^{2}$, we have $k>i$, $l>j$,
then again one can verify easily that
$Z (\ee )=Z_{+}(\ee )$. 
\end{Remarks}

We now turn to the classification of almost principal nilpotent pairs.

\begin{Definition}
Recall that a nilpotent pair $\ee$ is called {\bf almost principal}
if $\dim Z (\ee )=\rank \ggoth+1$.
\end{Definition}

It is shown in \cite{[P2]} that almost principal nilpotent pairs are
distinguished and wonderful. Therefore the number of $G$-orbits of
almost principal nilpotent pairs is finite.

\begin{Theorem}
Let $\ggoth$ be of type $A$ or $D$, then
there are no almost principal nilpotent pairs.
\end{Theorem}
\begin{proof}
For type A, this is a direct consequence of Theorem 2.9 and
Proposition 4.2.

So let us suppose that $\ggoth$ is of type D, and let
$(\Gamma^{1},\Gamma^{2},\Gamma^{3})$ be the corresponding
triple of centrally symmetric skew diagram. If they are all
non empty, then by Proposition 3.11, $\E_{+}(\Gamma^{1},\Gamma^{2})$ and 
$\E_{+}(\Gamma^{1},\Gamma^{3})$ are not empty where
$\E_{+}(\Gamma^{1},\Gamma^{i})=
\bigcup_{p,q\geq 0}\E_{p,q}(\Gamma^{1},\Gamma^{i})$.
So $\Card Z (\ee )\geq \rank \ggoth+2$.

If $\Gamma^{2}=\Gamma^{3}=\emptyset$, then since $\Gamma^{1}$ is
non integral, $\Card Z (\ee )\neq \rank \ggoth+1$ by Proposition 4.5.

Finally, suppose that $\Gamma^{1}=\emptyset$. Then
the condition $\Card Z (\ee )= \rank \ggoth+1$ would
imply that by 
$$
\Card \E (\Gamma^{2},\Gamma^{2})+\Card \E (\Gamma^{3},\Gamma^{3})
+\Card \E (\Gamma^{2},\Gamma^{3})=\rank \ggoth+1
$$
where $\E (\Gamma^{2},\Gamma^{3})=\bigcup_{(p,q)\in\Z^{2}} 
\E_{p,q}(\Gamma^{2},\Gamma^{3})$.
So by Corollary 6.1, we deduce that
either $\Gamma^{2}$ or $\Gamma^{3}$ does not satisfy (R).
Without loss of generality, we can suppose that $\Gamma^{2}$ does not
satisfy (R). Thus we have by Theorem 3.5 and Proposition 4.4 that
$$
\Card \E (\Gamma^{2},\Gamma^{2})=1+\Card \E_{+} (\Gamma^{2},\Gamma^{2}),
\Card\E (\Gamma^{3},\Gamma^{3})=\Card \E_{+} (\Gamma^{3},\Gamma^{3}),
$$
$$
\Card \E (\Gamma^{2},\Gamma^{3})=1.\,\,\,\,\,\,\,\, (*)
$$
Since $\Gamma^{2}$ is integral, this means that 
it has at least two distinct southwest corners $(i,j),(k,l)$.
Let $(r,s)$ be a southwest corner of $\Gamma^{3}$. Then
the pairs $\{(i,j)\},\{(r,s)\}$ and $\{(k,l)\},\{(r,s)\}$
are in $\E (\Gamma^{2},\Gamma^{3})$. Hence, 
$\Card \E (\Gamma^{2},\Gamma^{3})\geq 2$ which contradicts $(*)$.
\end{proof}

\begin{Remarks}
In \cite{[P2]}, it is shown that there are no almost principal nilpotent
pairs for $\ggoth=E_{6}$, $E_{7}$, $E_{8}$ or $F_{4}$ and it is
suggested that almost principal nilpotent pairs do not exist in the 
simply-laced case. It is now clear by Theorem 6.4
that:
\end{Remarks}
\begin{Corollary} 
Almost principal nilpotent pairs do not exist in the simply-laced 
case. 
\end{Corollary}

This is rather strange and it would be nice to have a more natural 
explanation.

\begin{Theorem}
Let $\ggoth$ be of type $B_{n}$. Then there is a one-to-one
correspondence between the set of conjugacy classes of almost principal 
nilpotent pairs and the set of pairs of centrally symmetric 
skew diagrams $(\Gamma^{1},\Gamma^{2})$ satisfying:

\begin{Itemize}
\Item[a)] $\Card\Gamma^{1}+\Card\Gamma^{2} = 2n+1$.
\Item[b)] $\Gamma^{1}$ is integral and $\Gamma^{2}$ is non integral.
\Item[c)] We have either
\ItemItem[i)] $\Gamma^{1}$ is near rectangular of type a) or b) and
$\Gamma^{2}=\emptyset$;
\ItemItem[ii)] or $\Gamma^{1}=\{ (0,0)\}$ and $\Gamma^{2}$ is rectangular. 
\end{Itemize}
\end{Theorem}
\begin{proof}
Let $\ee$ be an almost principal nilpotent pair, and let 
$(\Gamma^{1},\Gamma^{2})$ be the associated pair of skew diagrams 
as given in Theorem 1.5. We have 2 cases:

\begin{Itemize}
\Item[1.] $\Gamma^{1}$ or $\Gamma^{2}$ is empty;
\Item[2.] $\Gamma^{1}$ and $\Gamma^{2}$ are both non-empty.
\end{Itemize}
\vskip0.5em

In the first case, by Proposition 4.5, we have condition c)i).
Let us show that in the 
second case, we have c)ii). So let us assume that
$\Gamma^{1}$ and $\Gamma^{2}$ are both non-empty.

By Theorem 3.5 and Proposition 3.11, we must have 
$\Card \E (\Gamma^{1},\Gamma^{2})=1$ and both
$\Gamma^{1}$ and $\Gamma^{2}$ must be rectangular for otherwise
one of them will have two distinct southwest corners and one can
use the argument in the proof of Theorem 6.4 for type D to obtain
a contradiction.

Now if $C_{1}$ (resp. $C_{2}$) is a skew subdiagram of $\Gamma^{1}$
(resp. $\Gamma^{2}$) such that $C_{1}$ and $C_{2}$ are of the same
shape, then one sees easily from the fact that the $\Gamma^{i}$'s
are centrally symmetric that the pair $(C_{1},C_{2})$ is in 
$\E (\Gamma^{1},\Gamma^{2})$. It follows that if
$\Card \E (\Gamma^{1},\Gamma^{2})=1$, then there is only one
such pair. Now $\Gamma^{2}$ is non integral,
therefore $\Gamma^{2}$ contains the square $\{ \pm 1/2,\pm 1/2\}$, hence
$\Gamma^{1}$ and $\Gamma^{2}$ must satisfy condition c)ii).

Finally, given $(\Gamma^{1},\Gamma^{2})$ satisfying a), b) and
c). It is clear that the corresponding distinguished nilpotent
pair is almost principal.  
\end{proof}

\begin{Theorem}
Let $\ggoth$ be of type $C_{n}$. Then there is a one-to-one
correspondence between the set of conjugacy classes of almost principal 
nilpotent pairs and the set of pairs of centrally symmetric 
skew diagrams $(\Gamma^{1},\Gamma^{2})$ satisfying:

\begin{Itemize}
\Item[a)] $\Card\Gamma^{1}+\Card\Gamma^{2} = 2n$.
\Item[b)] $\Gamma^{1}\subset \Z^{2}+(0,1/2)$ and 
$\Gamma^{2}\subset \Z^{2}+(1/2,0)$.
\Item[c)] We have either
\ItemItem[i)] $\Gamma^{1}$ is near rectangular and
$\Gamma^{2}=\emptyset$ or vice versa;
\ItemItem[ii)] or $\Gamma^{1}$ and $\Gamma^{2}$ are both (non-empty) 
rectangular with $\Gamma^{1}\subset \{0\} \times (\Z + 1/2 )$ and  
$\Gamma^{2}\subset (\Z +1/2 )\times \{ 0 \}$.
\end{Itemize}
\end{Theorem}
\begin{proof}
The proof is similar to the one for type B. The only difference
is that $\Gamma^{1}$ and $\Gamma^{2}$ are semi-integral which 
gives condition c)ii). 
\end{proof}

\section{Centralizer of almost principal nilpotent pairs}

In \cite{[P2]}, almost principal nilpotent pairs satisfying condition c)i) 
(resp. condition c)ii)) of Theorem 6.7 or Theorem 6.8 are called of
$\Z$-type (resp. non-$\Z$-type). He proved that the centralizer of
an almost nilpotent pair of non-$\Z$-type is abelian and suggested
the same to be true for an almost nilpotent pair of $\Z$-type.
We shall prove in this section that this is indeed true. 

\begin{Proposition}
Let $\ee$ be a distinguished nilpotent pair in a classical simple
Lie algebra $\ggoth$. 
Then there exists an abelian subalgebra in $Z_{+}(\ee )_{int}$
of dimension $\rank \ggoth$.
\end{Proposition}
\begin{proof}
Let $\Gamma$ be a skew diagram (resp. centrally symmetric)
and $C\in \E_{p,q}(\Gamma )$ 
(resp. $(C,\widetilde{C}\in \E_{p,q}(\Gamma,\Gamma)$),
$C'\in \E_{p',q'}(\Gamma )$ where $p,q,p',q'\geq 0$
(resp. $(C',\widetilde{C}'\in \E_{p,q}(\Gamma,\Gamma)$). 

Let $(i,j)\in C$ be such that $(i,j)+(p,q)\in C'$ 
(resp. $\in \widetilde{C}'$), then by the 
definition of a skew subdiagram and the fact that $p,q\geq 0$, we have
$(i,j)\in C'$ (resp. $\in \widetilde{C}'$). 

Now let $\Gamma,\widetilde{\Gamma}$ be integral centrally symmetric
skew diagrams such that $\Gamma \cap \widetilde{\Gamma}=\{ (0,0\}$ and 
$\widetilde{\Gamma}\neq \{ (0,0)\}$. 
By the proof of Theorem 3.9, there exists a pair 
$(C,\widetilde{C})\in \E_{p,q}(\Gamma ,\widetilde{\Gamma})$, $p,q\geq 0$ 
such that $C=\{ (i,j)\}$ and $\widetilde{C}=\{ (k,l)\}$ are 
southwest corners of $\Gamma$ and $\widetilde{\Gamma}$ respectively. 
Thus neither $(i,j)+(p,q)$ or $(k,l)+(p,q)$ is contained in $C'$ or 
$\widetilde{C}'$ for any pair 
$(C',\widetilde{C}')\in \E_{+} (\Gamma ,\Gamma )\cup
\E_{+}(\widetilde{\Gamma},\widetilde{\Gamma})$.

Now one verifies easily from Theorem 2.9 that the subalgebra
spanned by the elements $x_{C}$ or $x_{C,C'}$ where
\vskip0.5em

Type A -- $C\in \E_{+}(\Gamma )\setminus \E_{0,0}(\Gamma )$,

Types B or C -- $(C,C')\in \E_{+}(\Gamma^{1},\Gamma^{1})\cup 
\E_{+}(\Gamma^{2},\Gamma^{2})$,

Type D -- $(C,C')\in \E_{+}(\Gamma^{1},\Gamma^{1})\cup 
\E_{+}(\Gamma^{2},\Gamma^{2}) \cup \E_{+}(\Gamma^{3},\Gamma^{3})$
and the pair in $\E_{+}(\Gamma^{2},\Gamma^{3})$ as given in the
proof of Theorem 3.9 which was described above,
\vskip0.5em

\noindent form an abelian subalgebra of $Z (\ee )$. 
By Theorem 3.2 and Theorem 3.5, its dimension
is $\rank \ggoth$.
\end{proof}

\begin{Theorem}
Let $\ee$ be an almost princiapl nilpotent pair in a classical
simple Lie algebra $\ggoth$. Then $Z_{+}(\ee )$ is abelian.
\end{Theorem}
\begin{proof}
By the classification (Theorems 6.4, 6.7 and 6.8) and Proposition
7.1, this is just a straightforward case by case computation.
\end{proof}

\end{document}